\input amstex
\documentstyle{amsppt}

\magnification=\magstep1 \NoRunningHeads
\loadbold
\input epsf

\topmatter
\title Generic nonsingular Poisson suspension is of type $III_1$
\endtitle

\author
Alexandre I. Danilenko, Zemer Kosloff and Emmanuel Roy
\endauthor

\address
B. Verkin Institute for Low Temperature Physics \& Engineering
of Ukrainian National Academy of Sciences,
47 Lenin Ave.,
 Kharkiv, 61164, UKRAINE
\endaddress
\email            alexandre.danilenko\@gmail.com
\endemail

\address
Einstein Institute of Mathematics, Hebrew University of Jerusalem, Givat Ram. Jerusalem, 9190401, ISRAEL
\endaddress
\email
zemer.kosloff\@mail.huji.ac.il
\endemail

\address
Laboratoire Analyse, G\'eom\'etrie et Applications, CNRS UMR 7539,
Universit\'e Paris
	13, Institut Galil\'ee, 99 avenue Jean-Baptiste Cl\'ement F93430 Villetaneuse, FRANCE
\endaddress

\email 
roy\@math.univ-paris13.fr
\endemail

\thanks
The research of Z.K. was partially supported by ISF grant No. 1570/17.
\endthanks

\abstract
It is shown  that for a dense $G_\delta$-subset of the subgroup of nonsingular transformations (of a standard infinite $\sigma$-finite  measure space) whose  Poisson suspensions are nonsingular, the corresponding Poisson suspensions are ergodic and of Krieger's type $III_1$.
\endabstract

\endtopmatter

\document

\head 0. Introduction
\endhead

In this paper we continue to study nonsingular Poisson suspensions for nonsingular transformations of infinite Lebesgue spaces $(X,\goth B,\mu)$ that we initiated in our previous work \cite{DaKoRo}.
We showed there that there is a maximal subset Aut$_2(X,\mu)$ of the group Aut$(X,\mu)$
of all nonsingular transformations $T$ of $(X,\goth B,\mu)$ for which the Poisson suspension $T_*$ is well defined as a nonsingular invertible transformation of the associated Poisson  probability  space $(X^*,\goth B^*,\mu^*)$.
Moreover,  Aut$_2(X,\mu)$ is a Polish group in an appropriate topology $d_2$ which is stronger 
than the usual weak topology \cite{DaKoRo}.
The suspension $T_*$ admits an equivalent  invariant probability measure if and only if $T$ admits an equivalent invariant measure $\nu$ such that $\sqrt{\frac{d\nu}{d\mu}}-1\in L^2(X,\mu)$. 
In this paper we consider the problem:
\roster
\item"" {\it Let $T_*$ do not admit an equivalent probability measure. Can $T_*$ be  ergodic? If yes,  what is  the Krieger type of $T_*$?}
\endroster
We answer affirmatively the first question and contribute partly to the second one.
Though we are unable so far  to construct  a concrete example of an ergodic conservative $T$ of type $III$ whose Poisson suspension is ergodic and of type $III$,
we instead  utilize the Baire category tools to prove a stronger ``existence'' result.

\proclaim{Theorem A (main result)} The subset of all $T\in\text{\rom{Aut}}_2(X,\mu)$ such that  $T$ is ergodic and of type $III_1$ and $T_*$ is ergodic and of type $III_1$ is a dense $G_\delta$ in $(\text{\rom{Aut}}_2(X,\mu),d_2)$.
\endproclaim
To prove this theorem we first  construct a  concrete example of a totally dissipative $T$ with $T_*$ being ergodic and of type $III_1$.
The   construction is  motivated by  a recent progress in the theory of nonsingular Bernoulli shifts achieved in  \cite{Ko1} and \cite{DaLe}, because if $T$ is totally dissipative then $T_*$ is always a nonsingular Bernoulli shift.
However the aforementioned papers deal only with the shifts on $\{0,1\}^\Bbb Z$ while we encounter indeed with the product spaces  $A^\Bbb Z$ with $A$  uncountable.
This situation is considerably more difficult and we are not sure that the techniques developed in \cite{Ko1} and \cite{DaLe} extends to it in the full generality.
However, we need only a very particular case  which, in turn, can be reduced further to the shift on $(\Bbb Z_+^\Bbb N,\bigotimes_{n\in\Bbb Z}\kappa_n)$, where  $(\kappa_n)_{n\in\Bbb Z}$ is a specially selected sequence of Poisson distributions.
Then we prove that  $T_*$ is of type $III_1$ by showing that the Maharam extension of $T_*$ is conservative and has nonsingular property $K$. 
For that we use essentially properties of Skellam distributions (see Appendix A) and 
L{\'e}vy's continuity theorem in addition to the  theory of nonsingular endomorphisms and measurable  orbit theory that were utilized in \cite{Ko1} and \cite{DaLe}.

Secondly, we prove that the conjugacy class of $T$ is dense in Aut$_2(X,\mu)$.
 Main Theorem follows from that and an additional fact that
the subset of ergodic type $III_1$ transformations
in Aut$_2(X,\mu)$  is a dense $G_\delta$ (in $d_2$)  \cite{DaKoRo}.

The group Aut$_1(X,\mu)$  and a homomorphism $\chi:\text{Aut}_1(X,\mu)\to\Bbb R$ were introduced in \cite{Ne}.
We showed in \cite{DaKoRo} that   $\text{Aut}_1(X,\mu)\subset \text{Aut}_2(X,\mu)$,
$\text{Aut}_1(X,\mu)$ is a Polish group in a topology $d_1$ which is stronger than $d_2$,
and $\chi$ is continuous.
If $\chi(T)\ne 0$ then $T$ is dissipative.
Thus, every ergodic transformation from $\text{Aut}_1(X,\mu)$ is contained in the proper closed subgroup Ker\,$\chi$.
The following statement is proved in the same way as Theorem~A.

\proclaim{Theorem B} The subset of all $T\in\text{\rom{Ker}}\,\chi$ such that  $T$ is ergodic and of type $III_1$ and $T_*$ is ergodic and of type $III_1$ is a dense $G_\delta$ in $(\text{\rom{Ker}}\,\chi,d_1)$.
\endproclaim

It is interesting to note that dynamical properties of $T_*$ are determined  not only by the dynamical properties of $T$ but also by the choice of $\mu$ inside its equivalence class.
Indeed, if $(X,\mu, T)$ is totally dissipative then we can find three different measures $\mu_1\sim\mu_2\sim\mu_3$  in the equivalence class of $\mu$ such that $(X^*,\mu_1^*,T_*)$ is a probability preserving Bernoulli shift, $(X^*,\mu_2^*,T_*)$ is  an ergodic system of type $III_1$ and $(X^*,\mu_3^*,T_*)$ is a totally dissipative system.
We ``refine'' further this phenomenon in a rather surprising way:
for each $T\in\text{Aut}_2(X,\mu)$ and $t>0$ we consider a dynamical system
$(X^*,\mu_t^*,T_*)$, where $\mu_t$ is the scaling of $\mu$ by $t$, i.e. $\mu_t(A):=t\mu(A)$
for each Borel subset $A\subset X$. 
Of course, $\frac{d\mu_t\circ T}{d\mu}=\frac{d\mu\circ T}{d\mu}$ for each $t>0$.
We then encounter with the following {\it phase transition} phenomenon.

\proclaim{Theorem C}
Let    $T\in\text{\rom{Aut}}_1(X,\mu)$.
If
there is  $\alpha>1$ such that $\alpha^{-1}<(T^n)'(x)<\alpha$ for each $n>0$ at a.e. $x\in X$
then there is $t_0\in[0,+\infty]$ such that
the Poisson suspension $(X^*,\mu_t^*,T_*)$ is conservative for each $t\in(0,t_0)$ and 
the Poisson suspension $(X^*,\mu_t^*,T_*)$ is totally dissipative for each $t\in (t_0,+\infty)$.
\endproclaim

The most interesting case is when the {\it bifurcation point} $t_0$ is {\it proper}, i.e.
$0<t_0<+\infty$.

\example{Example D} There is a totally dissipative $T\in\text{\rom{Aut}}_1(X,\mu)$ and 
$t_0\in(\frac16,4)$ such that
 the Poisson suspension $(X^*,(t\mu)^*,T_*)$ is weakly mixing of stable type $III_1$ if $0<t<t_0$ and
 the Poisson suspension $(X^*,(t\mu)^*,T_*)$ is totally dissipative for $t>t_0$.
\endexample

As a byproduct, we apply the techniques developed in this paper to study  conservativeness of the Poisson suspensions of general dissipative transformations. 
In particular, we show the following.

\proclaim{Theorem E} If $T\in\text{\rom{Aut}}_1(X,\mu)$ and $\chi(T)\ne 0$ then $T_*$ is totally dissipative.
\endproclaim

The outline of the paper is as follows.
The first two sections are of preliminary nature.
We present there some concepts and facts from  the theory of measured equivalence relations and their cocycles (\S1) and nonsingular endomorphisms and their extensions (\S2) to be used below in the paper.
In \S\,3 we study Poisson suspensions of transformations  defined on purely atomic measure spaces.
Since the purely atomic case was not considered in our previous paper on the Poisson suspensions, we first establish some basic results related to the nonsingularity and conservativeness of such suspensions independently of  \cite{DaKoRo} (see Propositions~3.1 and 3.2).
Then we prove some necessary (Proposition~3.3) and sufficient (Proposition~3.4) conditions for conservativeness of the Poisson suspensions in terms of the underlying dynamical system.
Theorem~3.5 plus Corollary~3.7 provide some conditions on a measure on $\Bbb Z$ under which  the Poisson suspension of the underlying unit translation on $\Bbb Z$ possesses the nonsingular property $K$. 
This is technically the most involved result of the paper.
Example~3.8 gives a concrete example of  a measure on $\Bbb Z$ satisfying those conditions. 
In \S4 we show how to pass from the purely atomic case studied in \S\,3 to the continuous case.
In particular, we construct a totally dissipative transformation on a nonatomic Lebesgue space whose Poisson suspension is weakly mixing and of type $III_1$ (see Theorem~4.2 and a remark just below it).
In \S\,5 we prove Theorems~A and B (see Theorem~5.3). \S\,6 is devoted to Poisson extensions of general dissipative transformations.
We prove there Theorem~E (see Theorem~6.1).
Some extension of Theorem~E to the more general case where $T\in\text{Aut}_2(X,\mu)$ is also discussed in that section (see Remark~6.2).
In \S\,7 we study the phase transitions for the conservativeness of Poisson suspensions when scaling the underlying intensity.
Theorem~C is proved there. 
 Example~D is also provided in \S\,7.
The paper has 
Appendix~A which  is devoted completely to the Skellam  distributions and their properties that we utilize in the proof of Theorem~3.5.

\head
1. Measured equivalence relations and their cocycles
\endhead

Let $(X,\goth B,\mu)$ be a standard $\sigma$-finite measure space. A Borel equivalence relation $\Cal R\subset X\times X$ is {\it countable} if for each $x\in X$, the $\Cal R$-equivalence class $\Cal R(x)$ is countable.
If $\Cal R$ is countable and  $A\in\goth B$ then the $\Cal R$-{\it saturation} $\Cal R(A):=\bigcup_{x\in A}\Cal R(x)$ of $A$ belongs to $\goth B$.
An $\Cal R$-saturated subset is also called {\it $\Cal R$-invariant}.
If for each  Borel subset $A$ of zero measure,  $\Cal R(A)$ is also of zero measure then $\Cal R$ is called $\mu$-{\it nonsingular}.
If the sub-$\sigma$-algebra $\{\Cal R(A)\mid A\in\goth B\}$ of $\Cal R$-invariant Borel subsets is trivial (mod 0) then $\Cal R$ is called {\it $\mu$-ergodic}.

From now on $\Cal R$ is countable and $\mu$-nonsingular.
The {\it full group} $[\Cal R]$ of $\Cal R$ consists of all nonsingular transformations $R$ of $(X,\mu)$ such that $Rx\in\Cal R(x)$ at a.e. $x\in X$.
Given a locally compact Polish group $G$, a Borel map $\alpha:\Cal R\to G$ is called a {\it cocycle} of $\Cal R$ if there is a Borel subset $N$ of zero measure such that
$$
\alpha(x,y)\alpha(y,z)=\alpha(x,z)\quad\text{for all  $x,y\in X\setminus N$ such that 
$
(x,y),(y,z)\in \Cal R$.}
$$
In the later sections of this paper we deal only with the cases where $G$ is either $\Bbb R$ or the multiplicative group $\Bbb R_+^*$.
A cocycle $\beta:\Cal R\to G$ is {\it cohomologous} to $\alpha$ if there is a Borel function
$\phi:X\to G$ and a Borel subset $N$ of zero measure such that
$$
\beta(x,y)=\phi(x)\alpha(x,y)\phi(y)^{-1}\quad\text{for all  $x,y\in X\setminus N$. }
$$
Fix a left Haar measure $\lambda_G$ on $G$.
The {\it $\alpha$-skew product} equivalence relation $\Cal R(\alpha)$ on the product space $(X\times G,\mu\times\lambda_G)$ is defined by:
$$
(x,g)\sim(y,h)\quad\text{if}\quad (x,y)\in\Cal R\quad\text{and}\quad h=\alpha(x,y)g.
$$
This equivalence relation is countable and $(\mu\times\lambda_G)$-nonsingular.
If $\Cal R(\alpha)$ is ergodic then $\alpha$ is called {\it ergodic}.
Of course, if $\alpha$ is ergodic then $\Cal R$ is ergodic.
If a cocycle $\beta$ is cohomologous to $\alpha$ and $\alpha$ is ergodic then $\beta$ is also ergodic.

We now isolate an important cocycle of $\Cal R$ with values in the multiplicative group $\Bbb R^*_+$.
It  is called the {\it  Radon-Nikodym cocycle} of $\Cal R$ and 
denoted by $\Delta_\Cal R$.
To define it,  we first fix a countable subgroup $\Gamma$ of Borel bijections of $X$ that generates $\Cal R$.
Such a group exists according to \cite{FeMo}.
Now we set
$$
\Delta_\Cal R(x,\gamma x):=\frac{d\mu\circ\gamma}{d\mu}(x),\quad x\in X, \gamma\in\Gamma.
$$
It follows from the chain rule for the Radon-Nikodym derivatives
 that $\Delta_\Cal R$ is a cocycle of $\Cal R$.
This cocycle is well defined, i.e. it does not depend on the choice of $\Gamma$ generating $\Cal R$\footnote{In this paper we do not distinguish between  objects (such as  subsets, maps, cocycles, etc.)
which agree almost everywhere.}.

\example{Example 1.1} Let $A$ be a countable set and let $\lambda_n$ be a non-degenerated distribution on $A$, i.e.  $\lambda_n(a)>0$ if $a\in A$ and $\sum_{a\in A}\lambda_n(a)=1$, for each $n\ge 1$.
We set $(X,\lambda):=(A^\Bbb N,\bigotimes_{n=1}^\infty\lambda_n)$.
Denote by $\Cal S$ the {\it tail equivalence relation} on $X$, i.e. two points $x=(x_n)_{n=1}^\infty$ and $y=(y_n)_{n=1}^\infty$ from $X$ are $\Cal S$-equivalent if there is $N>0$ such that
$x_n=y_n$ for all $n>N$.
Then $\Cal S$ is an ergodic  $\lambda$-nonsingular countable equivalence relation on $X$
and
$$
\Delta_\Cal S(x,y)=\prod_{n=1}^\infty\frac{\lambda_n(y_n)}{\lambda_n(x_n)}=\prod_{n=1}^N\frac{\lambda_n(y_n)}{\lambda_n(x_n)}
$$
for all $(x,y)\in\Cal S$.
\endexample

If we change $\mu$ with an equivalent $\sigma$-finite measure then the Radon-Nikodym cocycle of $\Cal R$ related to the new measure is cohomologous to the original $\Delta_\Cal R$.
The  $\Delta_\Cal R$-skew product equivalence relation $\Cal R(\Delta_\Cal R)$
is called the {\it Maharam extension} of $\Cal R$.
If it is ergodic, i.e. $\Delta_\Cal R$ is ergodic, then $\Cal R$ is said {\it to be of Krieger's type $III_1$.}
It follows from the aforementioned cohomology property of $\Delta_\Cal R$ that the property of $\Cal R$ to be of type $III_1$ does not change if we replace $\mu$ with an equivalent $\sigma$-finite measure.

Suppose now that $\Cal R$ is $\mu$-ergodic.
Given a cocycle  $\alpha$ of $\Cal R$ with values in an Abelian locally compact Polish group  $G$, we say that an element $g\in G$ is {\it an essential value of $\alpha$} if for each  neighborhood $U$ of $g$ in $G$ and each Borel subset $A\subset X$ of positive measure there are a Borel subset $B\subset A$ of positive measure and a Borel one-to-one map $\gamma:B\to A$ such that $(x,\gamma x)\in\Cal R$ and  $\alpha(x,\gamma x)\in U$ for each $x\in B$.
The set  $r(\alpha)$ of all essential values of $\alpha$ is a closed subgroup of $G$ (\cite{Sc}, \cite{FeMo}).
The cocycle $\alpha$ is ergodic if and only if $r(\alpha)=G$ (\cite{Sc}, \cite{FeMo}).

The following standard approximation lemma (see \cite{ChHaPr, Lemma~2.1}) is useful for computation of $r(\alpha)$.

\proclaim{Lemma 1.2} Let $\goth A\subset\goth B$ be a semiring such the corresponding ring is dense in 
$\goth B$.
Let $0<\delta<1$ and let $g\in G$.
If for each $A\in\goth A$ of positive measure and a neighborhood $U$ of $g$ there are a subset $B\subset A$ and a one-to-one Borel map $\gamma:B\to A$ such that $\mu(B)>\delta\mu(A)$,
$(x,\gamma x)\in\Cal R$, $\alpha(x,\gamma x)\in U$ and $\delta<\Delta_{\Cal R}(x,\gamma x)<\delta^{-1}$ for all $x\in B$ then $g$ is an essential value of $\alpha$.
\endproclaim

An ergodic invertible transformation $T$ of $(X,\goth B,\mu)$ is called {\it of Krieger's type $III_1$} if the $T$-orbit equivalence relation is of type $III_1$.
If, moreover, $T\times S$ is of type $III_1$ for each ergodic probability preserving transformation $S$ then $T$ is called {\it of stable type $III_1$}.

For more information on the measurable orbit theory and detailed proofs of the aforementioned facts we refer the reader to \cite{FeMo}, \cite{Sc}, \cite{DaSi}.

\head 2. Nonsingular endomorphisms, their extensions and associated equivalence relations\endhead
Let $(X,\goth B,\mu)$ be a standard $\sigma$-finite measure space.
A Borel map $T:X\to X$ is called a {\it $\mu$-nonsingular endomorphism} if $\mu\circ T^{-1}\sim\mu$.
Consider the following decreasing sequence  
$
\goth B\supset T^{-1}\goth B\supset T^{-2}\goth B\supset\cdots
$
 of sub-$\sigma$-algebras in $\goth B$.
If $\bigwedge_{n=1}^{+\infty} T^{-n}\goth B=\{\emptyset, X\}$ (mod 0) then $T$ is called {\it exact}.
We will consider only {\it aperiodic} endomorphisms, i.e. we assume that 
$$
\mu\left(\bigcup_{n>0}\{x\in X\mid T^nx=x\}\right)=0.
$$
Let $T$ be $\mu$-nonsingular endomorphism.
We  recall that a  measurable function $\omega:X\to\Bbb R_+^*$ is called {\it markovian } for
$(X,\mu,T)$  if
$$
\int f\circ T\omega\,d\mu=\int f\,d\mu\quad\text{for each $f\in L^1(X,\mu)$  \cite{SiTh1}.}
$$
Such a function may not be unique (see \cite{SiTh1}).
However if we assume that 
 the measure $\mu\circ T^{-1}$ is $\sigma$-finite then a standard verification shows that there exists a unique $T^{-1}\goth B$-measurable markovian function for $(X,\mu,T)$. 
 We call it the {\it  Radon-Nikodym derivative of $T$} and denote by $\rho_T$.
 One can check that
 $\rho_T=\frac{d\mu}{d\mu\circ T^{-1}} \circ T$.
 Of course, if $T$ is invertible then $\rho_T=\frac{d\mu\circ T}{d\mu}$.

Let $\kappa$ denote a measure on $\Bbb R^*_+$ equivalent to the Lebesgue measure and such that $\kappa(aB)=a^{-1}\kappa(B)$.
Given $T$ such that $\mu\circ T^{-1}$ is $\sigma$-finite,
 we can define a $\sigma$-finite measure preserving endomorphism $T_{\rho_T}$ of the product space $(X\times\Bbb R^*_+,\mu\times\kappa)$ by setting $T_{\rho_T}(x,t):=(Tx,\rho_T(x)t)$. It is called {\it the Maharam extension of $T$}.
If  $(Y,\goth Y,\nu)$ is a standard $\sigma$-finite measure space, $S$ is a transformation from $\text{Aut}(Y,\nu)$, $\pi:Y\to X$ is a Borel map such that 
$$
\text{$\mu\circ\pi^{-1}=\nu$, $\pi S =T\pi$, $\rho_S=\rho_T\circ \pi$ \ and \ $\goth Y=\bigvee_{n>0}S^n\pi^{-1}\goth B$}
$$
 then $S$ is called {\it the natural extension of $T$} \cite{SiTh2}. 
 The existence and uniqueness (up to a natural isomorphism) of the natural extension was 
 proved in  \cite{SiTh2}.
We will denote it by $\widetilde T$.
The natural extension of the Maharam extension of $T$ is canonically isomorphic to the
Maharam extension of the natural extension of $T$.
If $S\in\text{Aut}(Y,\nu)$ is the natural extension of an exact endomorphism then $S$ is called {\it a nonsingular $K$-automorphism}.
We recall that a nonsingular invertible transformation is called {\it weakly mixing} if the Cartesian product of it with every ergodic probability preserving transformation is ergodic.
A nonsingular $K$-automorphism is weakly mixing whenever it is conservative \cite{SiTh2}.

If for a.e. $x\in X$, the set $T^{-1}\{x\}$ is at most countable then $T$ is called {\it countable-to-one}.
From now on we will consider only endomorphisms which are countable-to-one.
Given such a $T$, we can associate an equivalence relation $\Cal S_T$ on $X$  by setting:
$$
(x,y)\in\Cal S_T\quad\text{if there is $n\ge 0$ such that $T^nx=T^ny$}.\tag 2-1
$$
Then $\Cal S_T$ is  countable and $\mu$-nonsingular.
It is ergodic if and only if  $T$ is exact \cite{Ha}. 
If~\thetag{2-1} holds, we set
$$
\alpha_{\rho_T}(x,y):=\rho_T(x)\cdots\rho_T(T^{n-1}x)\rho_T(T^{n-1}y)^{-1}\cdots\rho_T(y)^{-1}.
$$
Then  $\alpha_{\rho_T}$ is a well defined cocycle\footnote{The reader should not confuse $\alpha_{\rho_T}$ with $\Delta_{\Cal S_T}$.} of $\Cal S_T$ with values in $\Bbb R^*_+$.

\remark{Remark \rom{2.1}}
We note that $\alpha_{\rho_T}$ depends also on $\mu$, i.e., in fact, $\alpha_{\rho_T}=\alpha_{\rho_T,\mu}$.
If we replace $\mu$ with an equivalent measure $\lambda$ and the Radon-Nikodym derivative $\frac{d\lambda}{d\mu}$ is measurable with respect to the $\sigma$-algebra $T^{-1}\goth B$ then the cocycle $\alpha_{\rho_T,\lambda}$ is cohomologous to $\alpha_{\rho_T,\mu}$.
\endremark

Denote by $\Cal S_T(\alpha_{\rho_T})$ the $\alpha_{\rho_T}$-skew product extension of $\Cal S$.
Then $\Cal S_{T_{\rho_T}}=\Cal S_T(\alpha_{\rho_T})$.
From the aforementioned facts we deduce the following proposition (see \cite{DaLe} for details).

\proclaim{Proposition 2.2}  Let $T$ be countable-to-one nonsingular endomorphism of a probability space.
The Maharam extension of the natural extension of $T$ is a $K$-automorphism
if and only if the cocycle $\alpha_{\rho_T}$ of $\Cal S_T$ is ergodic. 
\endproclaim

We illustrate the aforementioned concepts with the following example.

\example{Example 2.3}
Let  $(X,\mu):=(A^\Bbb N,\bigotimes_{n=1}^\infty\mu_n)$, where $A$ and $(\mu_n)_{n=1}^\infty$ stand for the same objects as in Example~1.1.
Denote by $T:X\to X$ the one-sided Bernoulli shift on $X$.
Of course, it is countable-to-one.
By the Kakutani criterion \cite{Ka}, $T$ is $\mu$-nonsingular if and only if 
$$
\sum_{n=1}^\infty\sum_{a\in A}\Big(\sqrt{\mu_n(a)}-\sqrt{\mu_{n+1}(a)}\Big)^2<\infty.
$$
Moreover, if $T$ is nonsingular then we can compute the Radon-Nikodym derivative of $T$:
$$
\rho_T(x)=\prod_{n=1}^\infty\frac{\mu_{n}(x_{n+1})}{\mu_{n+1}(x_{n+1})}\qquad\text{at $\mu$-a.e. $x=(x_n)_{n\in\Bbb N}\in X$.}\tag2-2
$$
It is easy to see that $\Cal S_T$ is the tail equivalence relation on $X$.
Of course, $T$ is exact.
The natural extension of $T$ is the two-sided Bernoulli shift on the space $(A^\Bbb Z,\bigotimes_{n\in\Bbb Z}\widetilde\mu_n)$, where $\widetilde\mu_n:=\mu_n$  if $n>0$
and $\widetilde\mu_n:=\mu_1$ if $n\le 0$.
The corresponding projection $\pi:A^\Bbb Z\to A^\Bbb N$ is given by the formula $\pi((x_n)_{n\in\Bbb Z}):=(x_n)_{n\in\Bbb N}$.
\endexample

Given a  standard probability space $(Y,\goth Y,\nu)$ and a transformation $S\in\text{Aut}(Y,\nu)$, we denote by $U_S$ the associated (with $S$) Koopman unitary operator in $L^2(Y,\nu)$, i.e.
$U_Sf:=f\circ S\cdot \sqrt{S'}$ for all $f\in L^2(Y,\nu)$\footnote{For the sake of simplicity,  we use here and below the notation $S'$ for the Radon-Nikodym derivative of $S$, i.e. $S':=\frac{d\nu\circ S}{d\nu}$.}.
We recall that $S$ is said to be
\roster
\item"---" {\it conservative} if for each subset $B\in\goth Y$ of positive measure, there is $n>0$ such that  $\nu(S^{-n}B\cap B)>0$,
\item"---" {\it dissipative} if it is not conservative,
\item"---" {\it totally dissipative} if there is  $B\in\goth Y$ such that $\bigsqcup_{n\in\Bbb Z} T^nB=Y$ (mod 0).
\endroster
The following  lemma from \cite{Ko1} (generalizing the corresponding result for nonsingular Bernoulli shifts from \cite{VaWa}) will be used  in the subsequent sections.

\proclaim{Lemma 2.4} If  $\sum_{n=1}^\infty\langle U_S^n1,1\rangle<\infty$ 
then $S$ is totally dissipative.
\endproclaim

We will also need  a sufficient condition for conservativeness of nonsingular transformations.

\proclaim{Proposition 2.5} Let $S$ be a nonsingular  transformation on a standard probability space $(Y,\nu)$.
Assume that $\frac{d\nu}{d\nu\circ S^{-n}}\in L^2(Y,\nu)$ for each $n>0$ and write
$b(n):=\big\|\frac{d\nu}{d\nu\circ S^{-n}}\big\|_2^2$.
If there is a sequence $(a(n))_{n=1}^\infty$ of positive reals such that $\sum_{n=1}^\infty a(n)=\infty$ but $\sum_{n=1}^\infty a(n)^2b(n)<\infty$ then $S$ is conservative.
\endproclaim
\demo{Proof} We set $A_n:=\{y\in Y\mid \frac{d\nu\circ S^{-n}}{d\nu}(y)<a(n)\}$.
By Markov's inequality,
$$
\nu(A_n)=\nu\bigg(\bigg\{y\in Y\,\bigg|\,\bigg( \frac{d\nu}{d\nu\circ S^{-n}}(y)\bigg)^2>\frac 1{a(n)^2}\bigg\}\bigg)\le a(n)^2b(n).
$$
Since $\sum_{n=1}^\infty a(n)^2b(n)<\infty$, it follows from the Borel-Cantelli lemma that
for $\nu$-a.e. $y$, there is $N(y)>0$ such that for each $n>N(y)$,
$$
 \frac{d\nu\circ S^{-n}}{d\nu}(y)\ge a(n).
$$
Now the condition $\sum_{n=1}^\infty a(n)=\infty$ yields that 
$\sum_{n=1}^\infty\frac{d\nu\circ S^{-n}}{d\nu}(y)=\infty.$
Hence $S$ is conservative by Hopf's criterion \cite{Aa}.
\qed
\enddemo

We also recall  definition of an r.f.m.p. extension for nonsingular maps.
A nonsingular endomorphism $T$ of $(X,\goth B,\mu)$ is called a {\it relatively finite measure preserving} (r.f.m.p.) extension of a nonsingular endomorphism $S$ of $(Y,\goth Y,\nu)$ if there is a Borel map $\pi:X\to Y$ such that
$\nu=\mu\circ \pi^{-1}$, $\pi T=S\pi$ and $\rho_T=\rho_S\circ \pi$.
For instance, the natural extension of a nonsingular endomorphism is r.f.m.p.

For more information on nonsingular endomorphisms we refer to \cite{Aa}, \cite{Ha}, \cite{SiTh2} and \cite{DaLe}.

\head 3.  Poisson suspensions of type $III_1$ over a discrete base 
\endhead

Let $(X,\goth B,\mu)$ be a $\sigma$-finite measure space and $\mu(X)=\infty$.
We first recall some definitions, notation and facts from \cite{DaKoRo}.
Aut$(X,\mu)$ stands for the group of all $\mu$-nonsingular invertible transformations.
To simplify the notation, for each $T\in\text{Aut}(X,\mu)$, we will denote the Radon-Nikodym derivative $\frac{d\mu\circ T^{-1}}{d\mu}$ of $T$ by $T'$.
We also let
$$
\align
\text{Aut}_2(X,\mu)&:=\{T\in\text{Aut}(X,\mu)\mid \sqrt{T'}-1\in L^2(X,\mu)\}\quad\text{and}\\
\text{Aut}_1(X,\mu)&:=\{T\in\text{Aut}(X,\mu)\mid T'-1\in L^1(X,\mu)\}.
\endalign
$$
The two sets are Borel subgroups of $\text{Aut}(X,\mu)$ endowed with the weak topology  and $\text{Aut}_1(X,\mu)\subset\text{Aut}_2(X,\mu)$.
For $j=1,2$, we
 define a topology $d_j$ on $\text{Aut}_j(X,\mu)$ by saying that a sequence $(T_n)_{n=1}^\infty$ of transformations from $\text{Aut}_j(X,\mu)$ converges to a transformation $T$ from
$\text{Aut}_j(X,\mu)$ if $T_n\to T$ weakly and $\|(T_n')^{1/j}-(T')^{1/j}\|_j\to 0$ as $n\to\infty$. 
Then $\text{Aut}_j(X,\mu)$ endowed with $d_j$ is a Polish group.
There exists a continuous homomorphism $\chi:\text{Aut}_1(X,\mu)\to\Bbb R$ defined by
the formula
$$
\chi(T):=\int_X(T'-1)\,d\mu.
$$

From now on and  till the end of the section let $X=\Bbb Z$ and let $T$ denote the unit translation, i.e. $Tn:=n+1$ for all $n\in\Bbb Z$.
Given a measure $\mu$ on $X$, 
we set $a_n:=\mu(n)$ for each $n\in\Bbb Z$.
Then $T$ is $\mu$-nonsingular if and only if $\mu$ is non-degenerated, i.e. $a_n>0$ for each $n\in\Bbb Z$.
Moreover,  
$$
T'(n)=\frac{a_{n-1}}{a_n}\quad\text{ for all $n\in\Bbb Z$.}
$$
Of course, $\text{Aut}(X,\mu)$ is the group of all permutations of $\Bbb Z$.
We see that
$T\in\text{Aut}_2(X,\mu)$ if and only if 
$$
\sum_{n\in\Bbb Z}\left(\sqrt{\frac{a_{n-1}}{a_n}}-1\right)^2a_n=\sum_{n\in\Bbb Z}(\sqrt{a_{n-1}}-\sqrt{a_n})^2<\infty.\tag 3-1
$$
In a similar way, $T\in\text{Aut}_1(X,\mu)$ if and only if $\sum_{n\in\Bbb Z}|a_{n-1}-a_n|<\infty.$
If the latter inequality holds then there exist the two limits $a_{+\infty}:=\lim_{n\to+\infty}a_n$
and $a_{-\infty}:=\lim_{n\to-\infty}a_n$.
Moreover, it is easy to verify that 
$$
\chi(T)=a_{+\infty}-a_{-\infty}.
\tag3-2
$$

Let $(X^*,\mu^*,T_*)$ denote the Poisson suspension of the dynamical system $(X,\mu,T)$.
The mapping $\omega\mapsto(\omega(n))_{n\in\Bbb Z}$ is an isomorphism of $(X^*,\mu^*)$
onto the infinite product space $(\Bbb Z_+^\Bbb Z,\bigotimes_{n\in\Bbb Z}\kappa_n)$, where 
$\kappa_n$ is the Poisson distribution with parameter $a_n$ for each $n\in\Bbb Z$.
Moreover, this mapping conjugates $T_*$ with the shift on $\Bbb Z_+^\Bbb Z$.
For this reason, from now on we view $T_*$ as the shift $\Bbb Z_+^\Bbb Z$ defined by
$$
(T_*y)_n:=y_{n+1}\quad \text{for all $y=(y_n)_{n\in\Bbb Z}\in \Bbb Z_+^\Bbb Z.$}
$$
It was shown in \cite{DaKoRo} that in the case where $\mu$ is non-atomic,  $T_*$ is $\mu^*$-nonsingular if and only if $T\in\text{Aut}_2(X,\mu)$.
We now verify that the same holds also  in our (purely atomic) case.

\proclaim{Proposition 3.1} $T_*$ is $\mu^*$-nonsingular if and only if \thetag{3-1} is satisfied. In this case for a.e. $y=(y_n)_{n\in\Bbb Z}\in X^*$,
$$
(T_*)'(y)=\prod_{n\in\Bbb Z}\frac{\kappa_{n-1}(y_n)}{\kappa_n(y_n)}.
$$
\endproclaim

\demo{Proof}
We first recall definition of the 
 Hellinger distance $H(\lambda,\xi)$ between two probability measures
 $\lambda$, $\xi$ on a countable set $C$:
$$
H^2(\lambda,\xi):=1-\int_{C}\sqrt{\frac{d\lambda}{d\xi}}\,d\xi=\frac1{2}
\sum_{c\in C}\Big(\sqrt{\lambda(c)}-\sqrt{\xi(c)}\Big)^2.
$$
Given two Poisson distributions $\nu_a$ and $\nu_b$ on $\Bbb Z_+$ with parameters $a$ and $b$ respectively, then the Hellinger distance $H(\nu_a,\nu_b)$ between $\nu_a$ and $\nu_b$ satisfies
$$
H^2(\nu_a,\nu_b)=1-e^{-\frac 12(\sqrt a-\sqrt b)^2}.
$$ 
Therefore it follows from the Kakutani criterion \cite{Ka} that $T_*$ is $\mu^*$-nonsingular if and only if
$$
\infty>2\sum_{n\in\Bbb Z}H^2(\kappa_n,\kappa_{n+1})=2\sum_{n\in\Bbb Z}(1-e^{-\frac12(\sqrt{a_n}-\sqrt{a_{n+1}})^2}).\tag3-3
$$
The series   in the righthand side of \thetag{3-3} converges 
if and only if \thetag{3-1} holds. \qed
\enddemo

In  a similar way, the Kakutani criterion \cite{Ka} can be applied to prove the following claim.

\proclaim{Proposition 3.2} 
Let $\lambda$ be another non-degenerated measure on $\Bbb Z$.
Then $\lambda^*$ is equivalent to $\mu^*$
if and only if $\sum_{n\in\Bbb Z}\Big(\sqrt{\lambda(n)}-\sqrt{\mu(n)}\Big)^2<\infty$.
\endproclaim

Our next purpose is to investigate when $T_*$ is conservative or dissipative in terms of the original system $(X,\mu,T)$.
Let $\Cal L_{-}(\mu)$ denote the set of limit points of the sequence $(a_n)_{n<0}$ and
let $\Cal L_{+}(\mu)$ denote the set of limit points of the sequence $(a_n)_{n>0}$. 
If  $T\in\text{Aut}_2(X,\mu)$ then \thetag{3-1} implies  that there exist  reals $\beta_+\ge\alpha_+\ge 0$ and $\beta_-\ge\alpha_-\ge 0$  such that $\Cal L_{-}(\mu)=[\alpha_-,\beta_-]$ and $\Cal L_{+}(\mu)=[\alpha_+,\beta_+]$.

\proclaim{Proposition 3.3}\roster
\item"\rom{(i)}" If $T\in\text{\rom{Aut}}_2(X,\mu)$ and $\Cal L_{-}(\mu)\cap\Cal L_{+}(\mu)=\emptyset$ then
$T_*$ is totally dissipative.
\item"\rom{(ii)}" If $T\in\text{\rom{Aut}}_1(X,\mu)$ and $\chi(T)\ne 0$ then $T_*$ is totally dissipative.
\endroster
\endproclaim

\demo{Proof} Since (ii) follows from (i) and \thetag{3-2}, it suffices to prove (i).
Let $H(.,.)$ denote the Hellinger distance on the set of probability measures on $\Bbb Z_+$.
It follows from the condition of (i) that there is $\delta>0$ and $N>0$ such that $H(\kappa_i,\kappa_j)>\delta$ whenever $i<-N$ and $j>N$. 
This inequality and Proposition~3.1 yield that if $n>3N$ then
$$
\align
\langle U_{T_*}^n1,1\rangle &=\prod_{n\in\Bbb Z}\int_{\Bbb Z_+}\sqrt{\frac{d\kappa_{k-n}}{d\kappa_k}(y_k)}\,d\kappa_k(y_k)\\
&=
\prod_{n\in\Bbb Z}(1-H^2(\kappa_{k-n},\kappa_k))\\
&\le
\prod_{\frac n3<k<\frac{2n}3}(1-H^2(\kappa_{k-n},\kappa_k))\\
&<(1-\delta^2)^{n/3}.
\endalign
$$
Hence $\sum_{n=1}^\infty\langle U_{T_*}^n1,1\rangle<\sum_{n=1}^\infty(1-\delta^2)^{n/3}<\infty$. 
 It remains to use Lemma~2.4.
\qed
\enddemo

We will also need a sufficient condition for conservativeness of $T^*$.

\proclaim{Proposition 3.4} Let $T\in\text{\rom{Aut}}_1(X,\mu)$, $\chi(T)=0$ and 
 $\Big(\frac{d\mu}{d\mu\circ T^{-n}}\Big)^2-1\in L^1(X,\mu)$  for each $n\in\Bbb N$.
 If there is a sequence $(b_n)_{n=1}^\infty$ of positive reals such that $\sum_{n=1}^\infty b_n=\infty$ but 
$$
\sum_{n=1}^\infty b_n^2e^{\int_X \big(\big(\frac{d\mu}{d\mu\circ T^{-n}}\big)^2-1\big)d\mu}<\infty
$$ 
then $T_*$ is conservative.
\endproclaim

\demo{Proof}
Since $\chi(T)=0$, it follows that $\sum_{k\in\Bbb Z}(a_k-a_{k-n})=0$ for each $n>0$.
Therefore
$$
\align
\log\bigg(\frac{d\mu^*}{d\mu^*\circ T_*^{-n}}(y)\bigg)^2&=-
{\sum_{k\in\Bbb Z}2\log\frac{\kappa_{k-n}(y_k)}{\kappa_k(y_k)}}\\
&=
-2\sum_{k\in\Bbb Z}\left(y_k\log\frac{a_{k-n}}{a_k}+a_k-a_{k-n}\right)\\
&=
-2\sum_{k\in\Bbb Z}y_k\log\frac{a_{k-n}}{a_k}=
{-\big\langle y,2\log\frac{d\mu\circ T^{-n}}{d\mu}\big\rangle}
\endalign
$$
at a.e. $y\in X^*$.
Here we consider $y$ as a measure on $X$ and use the notation $\langle y,f\rangle$ for the integral of a function $f\in L^1(X,y)$ with respect to $y$.
Utilizing the above formula and  the Laplace transform  we obtain  that
$$
\align
\bigg\|\frac{d\mu^*}{d\mu^*\circ T_*^{-n}}\bigg\|_2^2&=
\int_{X^*}e^{-\big\langle y,2\log\frac{d\mu\circ T^{-n}}{d\mu}\big\rangle}d\mu^*(y)\\
&=
e^{\int_X\Big(e^{-2\log\frac{d\mu\circ T^{-n}}{d\mu}}-1\Big)d\mu}\\
&=
e^{\int_X\Big(\big(\frac{d\mu}{d\mu\circ T^{-n}}\big)^2-1\Big)d\mu}.
\endalign
$$
It remains to apply Proposition~2.5.
\qed
\enddemo

From now on we will assume that $a_n:=ae^{\epsilon_n}$ for some  $a>0$ and a sequence $(\epsilon_n)_{n\in\Bbb Z}$ of reals such that 
$$
 \text{$\epsilon_{n}=0$ if $n\le 1$,  $\lim_{n\to+\infty}\epsilon_n= 0$,
 $\sum_{n=1}^\infty\epsilon^2_n=\infty$ but $\sum_{n=1}^\infty\epsilon^4_n<+\infty$.}\tag 3-4
 $$
Let $R$ denote the one-sided Bernoulli shift on the  space $((\Bbb Z_+)^\Bbb N,\bigotimes_{n=1}^\infty\kappa_n)$,  and let $\pi$ denote the
 projection 
 $$
 (\Bbb Z_+)^\Bbb Z\ni(x_n)_{n\in\Bbb Z}\mapsto(x_n)_{n>0}\in (\Bbb Z_+)^\Bbb N.
 $$
 Then $R\pi=\pi T_*$.
In view of Example~2.3,  $T_*$ is the natural extension of $R$.

\proclaim{Theorem 3.5}  If \thetag{3-4} holds then 
the cocycle $\log\alpha_{\rho_R}$ of the equivalence relation $\Cal S_R$ associated with $R$ is ergodic.
\endproclaim

\demo{Proof}
According to Example~2.3,
 $\Cal S_R$ is the tail equivalence relation on $(\Bbb Z_+)^\Bbb N$.
 By~\thetag{2-2},
 $\rho_R(x)=\prod_{j=2}^\infty\frac{\kappa_{j-1}(x_j)}{\kappa_j(x_j)}$ for each $x\in X$.
Hence for each $n>0$,
$$
\rho_R(x)\cdots\rho_R(R^{n-1}x)=\prod_{j=2}^{n}\frac{\kappa_1(x_j)}{\kappa_j(x_j)}\prod_{j>n}
\frac{\kappa_{j-n}(x_j)}{\kappa_j(x_j)}.
$$
Suppose now that $R^nx=R^ny$.
Then $x_j=y_j$ for all $j>n$.
Therefore
$$
\alpha_{\rho_R}(x,y):=\frac{\rho_R(x)\cdots\rho_R(R^{n-1}x)
}{\rho_R(y)\cdots\rho_R(R^{n-1}y)}=
\prod_{j=2}^{n}\frac{\kappa_{j}(y_j)}{\kappa_{j}(x_j)}
\prod_{j=2}^{n}\frac{\kappa_{1}(x_j)}{\kappa_{1}(y_j)}=\prod_{j=2}^n\bigg(\frac{a_j}{a_1}\bigg)^{y_j-x_j}.
$$
Thus we obtain that for each $(x,y)\in \Cal S_R$,
$$
\log \alpha_{\rho_R}(x,y)=\sum_{j>1}(y_j-x_j)\epsilon_j.
$$
Let $\Delta_{\Cal S_R}$ denote the Radon-Nikodym cocycle of $\Cal S_R$.
Utilizing Example~1.1 and the definition of $\kappa_j$ we obtain that 
$$
\log\Delta_{\Cal S_R}(x,y)=\log \bigg(\prod_{j=1}^{n}\frac{\kappa_{j}(y_j)}{\kappa_{j}(x_j)}\bigg)=\sum_{j>0}\bigg((y_j-x_j)\log a_j-\log\frac{y_j!}{x_j!}\bigg).
$$
Now we fix an infinite subset $J\subset\{2,3,4,\dots\}$ with $\sum_{j\in J}\epsilon_j^2<\infty$.
Then it follows from~\thetag{3-4} that $\sum_{j\not\in J}\epsilon_j^2=\infty$.
Hence we can
 change $(\epsilon_n)_{n\in\Bbb Z}$ with another sequence  coinciding with  the original one if $n\not\in J$ and such that 
 $\epsilon_n=0$ for each $n\in J$.
Of course, the new sequence
 satisfies \thetag{3-4} and
  determines a new measure $\widehat\mu$ on $X$.
 Proposition~3.2 yields that $\widehat\mu^*\sim\mu^*$.
 Moreover, the Radon-Nikodym derivative $x\mapsto\frac{d\widehat\mu^*}{d\mu^*}(x)$ does not depend
 on the first coordinate of $x=(x_j)_{j=1}^\infty$ because $1\not\in J$.
Therefore according to Remark~2.1, the cocycle $\alpha_{\rho_R,\widehat\mu^*}$
 is cohomologous to $\alpha_{\rho_R,\mu^*}$.
 Hence $\alpha_{\rho_R,\widehat\mu^*}$ is ergodic if and only if $\alpha_{\rho_R,\mu^*}$
 is ergodic.
 In the sequel, by $\alpha_{\rho_R}$ we mean $\alpha_{\rho_R,\widehat\mu^*}$.
Thus, 
 we may assume without loss of generality that  the triplet $(X^*,\mu^*,\Cal S_R)$
is isomorphic to  $(Z,\eta,\Cal T)$, where
$Z=(\Bbb Z_+\times\Bbb Z_+)^\Bbb N$, $\eta=\bigotimes_{j=1}^\infty(\kappa_j\otimes\kappa_0)$
and $\Cal T$ is the tail equivalence relation on $Z$\footnote{The corresponding isomorphism is given by the map $X^*\ni x\mapsto((x_j)_{j\not\in J},(x_j)_{j\in J})\in (\Bbb Z_+\times\Bbb Z_+)^\Bbb N$.}.
Let
 $z=(x_j,y_j)_{j>0}$ and $z'=(x_j',y_j')_{j>0}$ are two $\Cal T$-equivalent points with $x_j,y_j, x_j',y_j'\in\Bbb Z_+$.
Computing $\log \alpha_{\rho_R}$ and $\log\Delta_{\Cal S_R}$ in the ``new coordinates'' we obtain that
$$
\align
\log \alpha_{\rho_R}(z,z')&=\sum_{j>0}(x_j'-x_j)\epsilon_j\quad\text{and}\\
\log\Delta_{\Cal S_R}(z,z')&=\sum_{j>0}\bigg((x_j'-x_j)\log a_j+(y_j'-y_j)\log a_0-\log\frac{x_j'!}{x_j!}-\log\frac{y_j'!}{y_j!}\bigg).
\endalign
$$
Denote by $\tau$ the flip on $\Bbb Z_+\times\Bbb Z_+$, i.e. $\tau(x,y)=(y,x)$.
Define  a transformation $\tau_n$ of $Z$ by setting $\tau_n:=\underbrace{\tau\times\cdots\times \tau}_{n\text{ times}}\times 
\text{\,Id}\times\text{Id}\cdots$.
Of course, $\tau_n\in[\Cal T]$ and 
$$
\log \alpha_{\rho_R}(z,\tau_nz)=\sum_{j=1}^n(y_j-x_j)\epsilon_j=\log\Delta_{\Cal S_R}(z,\tau_nz)
\tag3-5
$$
at a.e. $z\in Z$ for each $n>0$.

{\it Claim 1.} For each $p>0$,
$$
\lim_{n\to\infty}\eta(\{z\in Z\mid \log \alpha_{\rho_R}(z,\tau_nz)>-p\})=0.
$$

To prove this claim we first define mappings $X_j:Z\to\Bbb R$ by setting
$$
X_j(z):=(y_j-x_j)\epsilon_j
$$
for each $z=(x_m,y_m)_{m>0}$ and $j\in\Bbb N$.
Then the following are satisfied:
\roster
\item"---"
$X_1,X_2,\dots$ is a sequence of independent random  variables,
\item"---"
 for each $j>1$, the measure $\eta\circ (\frac 1{\epsilon_j}X_j)^{-1}$ is the Skellam
distribution $\chi_{a_0,a_j}$ with parameters $(a_0,a_j)$ (see Appendix~A) and
\item"---" for each $n>0$ and a.e. $z\in Z$,
$$
\log \alpha_{\rho_R}(z,\tau_nz)=\sum_{j=1}^nX_j(z). 
\tag{3-6}
$$
\endroster
We now have (see Appendix A): 
$$
\align
E(X_j)&=\epsilon_j(a_0-a_j)=a\epsilon_j(1-e^{\epsilon_j}),\\
 \sigma^2(X_j)&=\epsilon_j^2(a_0+a_j)=a\epsilon_j^2(1+e^{\epsilon_j})
 \endalign
 $$
  and for each $t\in\Bbb R$,
$$
\aligned
\phi_{X_j-E(X_j)}(t)&:=E(e^{it(X_j-E(X_j))})\\
&=e^{-(a+ae^{\epsilon_j})+ae^{it\epsilon_j}+ae^{\epsilon_j}e^{-it\epsilon_j}-itE(X_j)}\\
&=e^{a((e^{it\epsilon_j}-1+e^{-it\epsilon_j}-1)+(e^{\epsilon_j}-1)(e^{-it\epsilon_j}-1+it\epsilon_j))}\\
&=e^{-4a\sin^2\frac{t\epsilon_j}{2}+a(e^{\epsilon_j}-1)(e^{-it\epsilon_j}-1+it\epsilon_j)}.
\endaligned
\tag 3-7
$$
Let $\beta_n:=\big(\sum_{j=1}^n\epsilon_j^2\big)^{-\frac12}$.
Utilizing  \thetag{3-7}, we now compute the characteristic function of the random variable $Y_n:=\beta_n\sum_{j=1}^n(X_j-E(X_j))$ at a point $t\in\Bbb R$:
$$
\aligned
\phi_{Y_n}(t)&=\prod_{j=1}^n\phi_{X_j-E(X_j)}(\beta_nt)\\
&=
e^{-4a\sum_{j=1}^n\sin^2\frac{t\epsilon_j\beta_n}2+a\sum_{j=1}^n(e^{\epsilon_j}-1)(e^{-it\beta_n\epsilon_j}-1+it\epsilon_j\beta_n))}.
\endaligned
\tag 3-8
$$
Since the sequence $(\epsilon_n)_{n=1}^\infty$ is bounded and $\sum_{j=1}^\infty\epsilon_j^2=\infty$, it follows that 
$$
\beta_n\max_{1\le j\le n}|\epsilon_j|\to 0\quad\text{ as $n\to\infty$.}
$$
Therefore there is $C=C(t)>0$ such that
$$
\align
\bigg|\sum_{j=1}^n\sin^2\frac{t\epsilon_j\beta_n}2-\sum_{j=1}^n\bigg(\frac{t\epsilon_j\beta_n}2\bigg)^2\bigg| &\le C\beta_n^4\sum_{j=1}^n\epsilon_j^4\quad\text{and}
\tag 3-9\\
\bigg|\sum_{j=1}^n(e^{\epsilon_j}-1)(e^{-it\beta_n\epsilon_j}-1+it\epsilon_j\beta_n))\bigg|
&\le C\beta_n^2\sum_{j=1}^n |\epsilon_j|^3.
\tag 3-10
\endalign
$$
For each $\epsilon>0$, there is $N>0$ such that $|\epsilon_n|<\epsilon$ whenever $n>N$.
Therefore 
$$
\align
&\beta_n^4\sum_{j=1}^n\epsilon_j^4\le
\beta_n^4\sum_{j=1}^N\epsilon_j^4+
\frac{\epsilon^2\sum_{j=N+1}^n\epsilon_j^2}{(\sum_{j=1}^n\epsilon_j^2)^2}<2\epsilon
\ \text{ and }\\
&\beta_n^2\sum_{j=1}^n|\epsilon_j|^3\le
 \beta_n^2\sum_{j=1}^N|\epsilon_j|^3 +\frac{\epsilon\sum_{j=N+1}^n\epsilon_j^2}{\sum_{j=1}^n\epsilon_j^2}<2\epsilon
\endalign
$$
if $n$ is large enough.
It follows that  the lefthand sides in  \thetag{3-9} and \thetag{3-10} go to $0$ as $n\to\infty$.
Hence \thetag{3-9} yields  that
$\lim_{n\to\infty}\sum_{j=1}^n\sin^2\frac{t\epsilon_j\beta_n}2=\frac{t^2}4$ for each $t\in\Bbb R$.
From this, \thetag{3-8} and \thetag{3-10},
 we  deduce that
$$
\lim_{n\to\infty}\phi_{\beta_n\sum_{j=1}^n(X_j-E(X_j))}(t)=e^{-at^2}\quad\text{for each $t\in\Bbb R$.}
$$
However, the map $\Bbb R\ni t\mapsto e^{-at^2}$ is the characteristic function of 
a normal distribution.
Hence by   L{\'e}vy's continuity theorem,  $(\beta_n\sum_{j=1}^n(X_j-E(X_j)))_{n=1}^\infty$  converges in distribution to a Gaussian random variable. 
Since $\sum_{j=1}^n\epsilon_j^2\to \infty$, it follows that
$$
\lim_{n\to\infty} \beta_n\sum_{j=1}^nE(X_j)=\lim_{n\to\infty}\frac{a\sum_{j=1}^n\epsilon_j(1-e^{\epsilon_j})}{\sqrt{\sum_{j=1}^n\epsilon_j^2}}=-\infty.
$$
Therefore $\beta_n\sum_{j=1}^nX_j\to-\infty$ in distribution.
This, in turn, yields that $\sum_{j=1}^nX_j\to-\infty$ in distribution.
Thus Claim~1 is proved.

{\it Claim 2.} For $ \eta$-a.e. $z\in Z$ and each $\epsilon>0$, there is $N$ such that 
$|X_n(z)|<\epsilon$ for each $n>N$.
To prove this claim we let
 $$
 B_n:=\{z=(x_j,y_j)_{j\in\Bbb N}\in Z\mid |X_n(z)|>\sqrt{|\epsilon_n|}\}=\{z\in Z\mid |y_n-x_n|>|\epsilon_n|^{-1/2}\}.
 $$
Then $\eta(B_n)=\chi_{a_0,a_n}(\{k\in\Bbb Z\mid |k|>|\epsilon_n|^{-1/2}\})$.
According to~\thetag{A-1}, there is $L>0$ such that if $|\epsilon_n|^{-1/2}>L$ then
$$
\sum_{|k|>|\epsilon_n|^{-1/2}}\chi_{a_0,a_n}(k)<\epsilon_n^4.
$$
Hence $\eta(B_n)<\epsilon_n^4$.
Therefore $\sum_{n=1}^\infty\eta(B_n)<\sum_{n=1}^\infty\epsilon_n^4<\infty$ according to \thetag{3-4}.
Claim~2 follows from this via  the Borel-Cantelli lemma.

{\it Claim 3.} Each $r<-1$ is an essential value of $\log\alpha_{\rho_R}$.

We first recall  that given $d>0$ and elements $c_1,\dots,c_d\in \Bbb Z_+\times\Bbb Z_+$,
the {\it cylinder} $[c_1,\dots,c_d]_1^d$ is the set
$\{z=(z_n)_{n\in\Bbb Z}\in Z\mid z_j=c_j, 1\le j\le d\}$.
To prove the claim, take $\epsilon>0$, a positive integer $k$ and a cylinder 
$C=[(a_1,b_1),\dots, (a_k,b_k)]_1^k\subset Z$,
 where $a_1,b_1,\dots,a_k,b_k\in\Bbb Z_+$.
It follows from Claim~1 and Claim~2 that there are $N>M>k$ and a subset $I\subset (\Bbb Z_+\times\Bbb Z_+)^{N-k}$ such that the subset
$$
A:=\bigsqcup_{(z_{k+1},\dots,z_N)\in I}[(a_1,b_1),\dots,(a_k,b_k),z_{k+1},\dots,z_N]_1^N\subset C
$$
satisfies the following three  conditions: 
\roster
\item"---"
$\eta(A)>0.5\eta(C)$, 
\item"---"$
\max_{z\in A}|X_j(z)|<\epsilon$ whenever $N\ge j>M$ and 
\item"---"$
\max_{z\in A}\sum_{j=M+1}^NX_j(z)<-r$.
 \endroster
 For $z\in A$, let $l(z)$ be the smallest number $l>M$ such that  $\sum_{j=M+1}^lX_j(z)<-r$.
 Then $l(z)\le N$ and
 $$
 \Bigg|\sum_{j=M+1}^{l(z)}X_j(z)-r\Bigg|\le |X_{l(z)}(z)|<\epsilon.
 \tag 3-11
 $$
 We now set 
 $$
 \psi(z):=((a_1,b_1),\dots,(a_k,b_k),z_{k+1},\dots,z_M,\tau(z_{M+1}),\dots, \tau(z_{l(z)}),z_{l(z)+1},\dots)
 $$
 for each $z\in A$.
 Then $\psi(z)\in C$, $(z,\psi(z))\in\Cal T$ and 
 $$
 |\log\alpha_{\rho_R}(z,\psi(z))-r|<\epsilon\tag3-12
 $$
 in view of \thetag{3-11} and \thetag{3-6}.
 We now show that the map $\psi:A\ni z\mapsto\psi(z)\in C$ is one-to-one.
 Suppose that $\psi(z)=\psi(z')$ for some $z=(z_j)_{j\in\Bbb Z}\in A$ and 
 $z'=(z_j')_{j\in\Bbb Z}\in A$.
 If $l(z)=l(z')$ then obviously $z=z'$.
 Therefore suppose that $l(z)>l(z')$.
 Then the equality $\psi(z)=\psi(z')$ implies that $z_j=z_j'$ if $1\le j\le M$
 and $\tau(z_j)=\tau(z_j')$ if $M<j\le l(z')$.
 Hence $z_j=z_j'$ if $1\le j\le l(z')$.
 Therefore $-r>\sum_{j=M+1}^{l(z')}X_j(z')=\sum_{j=M+1}^{l(z')}X_j(z)$.
 This yields that $l(z)\le l(z')$, a contradiction.
 Thus, $\psi$ is one-to-one.
It follows from  and \thetag{3-5}  and \thetag{3-12} that
$-r-\epsilon<\log\Delta_{\Cal S_R}(z,\tau_nz)<-r+\epsilon$ for all $z\in A$.
Lemma~1.2 implies now that
 $r$ is an essential value of $\log\alpha_{\rho_R}$.
Thus, Claim~3 is proved.

Since the essential range of $\log\alpha_{\rho_R}$ is a closed subgroup of $\Bbb R$
containing an arbitrary real less than $-1$ (in view of Claim~3), it follows that this subgroup is $\Bbb R$. 
Hence
$\log\alpha_{\rho_R}$ is ergodic.
\qed
\enddemo

\remark{Remark 3.6} While proving Theorem~3.5 we also showed as a byproduct  that the cocycle $\Delta_{\Cal S_R}$ is also ergodic, i.e. $\Cal S_R$ is of type $III_1$.
In fact, we proved a stronger result.
Let $\Cal E$ denote the orbit equivalence relation on $(\Bbb Z_+)^{\Bbb N}$ generated by the group of finite permutations of the coordinates.
Then $\Cal E$ is a proper  subrelation of $\Cal S_R$.
It follows from the proof of Theorem~3.5 that $\Cal E$ is of type $III_1$ and 
the restriction of the cocycle $\log\alpha_{\rho_R}$ to $\Cal E$ is ergodic.
\endremark

The next corollary is the main result of this section. 

\proclaim{Corollary 3.7}  If \thetag{3-4} holds and $T_*$ is conservative then
 the Maharam extension of $T_*$  is a weakly mixing $K$-automorphism. 
 In particular, $T_*$ is  weakly mixing and of stable type $III_1$.
\endproclaim

\demo{Proof} 
 It follows from Theorem~3.5 that $\alpha_{\rho_R}$ is ergodic.
 Since $T_*$ is the natural extension of $R$,
 the Maharam extension  $(T_*)_{\rho_{T_*}}$ of $T_*$ is a nonsingular $K$-automorphism
 by Proposition~2.2.
Since $T_*$ is conservative, $(T_*)_{\rho_{T_*}}$ is also conservative according to the Maharam theorem (see \cite{Aa}, \cite{Sc}).
It follows that $(T_*)_{\rho_{T_*}}$ is weakly mixing.
This implies that $T_*$ is weakly mixing  and of type $III_1$. \qed
\enddemo

We now provide  concrete examples of conservative $T_*$ such that \thetag{3-4} holds.

\example{Example 3.8} Let $\epsilon_n:=0$ if $n\le 1$ and $\epsilon_n:=-n^{-1/2}$ if $n>1$.
Then \thetag{3-4} holds.
We claim that  if $0<a<\frac 1{6}$ then $T_*$ is conservative.
For that we will utilize Proposition~3.4.
We first note that
$$
\align
\int_X\left(\left(\frac{d\mu}{d\mu\circ T^{-n}}\right)^2-1\right)d\mu&=
a\sum_{k\in\Bbb Z}\Big(e^{2(\epsilon_k-\epsilon_{k-n})}-1\Big)e^{\epsilon_k}\\
&=
a\sum_{k\in\Bbb Z}\Big(e^{3\epsilon_k-2\epsilon_{k-n}}-e^{\epsilon_k}\Big)
\endalign
$$
and $|3\epsilon_k-2\epsilon_{k-n}|\le3|\epsilon_k|$.
Therefore it follows from the Taylor expansion of the exponential function  
that
$$
e^{3\epsilon_k-2\epsilon_{k-n}}-e^{\epsilon_k}=2\epsilon_k-2\epsilon_{k-n}+\frac12(3\epsilon_k-2\epsilon_{k-n})^2-\frac12\epsilon_k^2+O(\epsilon_k^3)
$$
as $k\to\infty$.
Since  $\sum_{k\in\Bbb Z}|\epsilon_k|^3<\infty$ and $\sum_{k\in\Bbb Z}(\epsilon_k-\epsilon_{k-n})=0$,
$$
\align
\sum_{k\in\Bbb Z}(e^{3\epsilon_k-2\epsilon_{k-n}}-e^{\epsilon_k})
&=\frac12\sum_{k\in\Bbb Z}\Big((3\epsilon_k-2\epsilon_{k-n})^2-\epsilon_k^2\Big)+O(1)\\
&=2\sum_{k\in\Bbb Z}\Big(  2\epsilon_k^2+\epsilon_{k-n}^2-3\epsilon_k\epsilon_{k-n}\Big) +O(1)\\
&=4\sum_{k=2}^{n+1}\frac1k+2\sum_{k>n+1}\bigg( \frac 1{k-n}- \frac 1{k}
+\frac 3k-\frac 3{\sqrt{k(k-n)}}\bigg) +O(1)\\
&=6\sum_{k=2}^{n+1}\frac1k+6\sum_{k>n+1} \bigg(\frac 1k-\frac 1{\sqrt{k(k-n)}}\bigg) +O(1)\\
&=6\log n+6\int_{n+1}^{+\infty}\bigg(\frac1t-\frac1{\sqrt{t(t-n)}}\bigg)\,dt+O(1)\\
&=6\log n+6\Big(\log t-\log\Big(t-\frac n2+\sqrt{t(t-n)}\Big)\Big)\bigg |_{n+1}^{+\infty}  +O(1)\\
&=6\log n+O(1)
\endalign
$$
as $n\to+\infty$.
Thus, we have shown that
$$
\int_X\left(\left(\frac{d\mu}{d\mu\circ T^{-n}}\right)^2-1\right)d\mu= 6a\log n+O(1).
$$
It follows that there is a constant $C>0$ such that for each $n>0$,
$$
\int_X\left(\left(\frac{d\mu}{d\mu\circ T^{-n}}\right)^2-1\right)d\mu<6a\log n+C.\tag3-13
$$
Choose  a real $\beta$ such that $1\ge\beta>\frac 12+3a$.
It exists because $a<\frac 1{6}$.
Now let $b_n:=n^{-\beta}$. 
Then $\sum_{n=1}^\infty b_n=\infty$ but in view of \thetag{3-13},
$$
\sum_{n=1}^\infty b_n^2e^{\int_X\left(\left(\frac{d\mu}{d\mu\circ T^{-n}}\right)^2-1\right)d\mu}
\le e^{C} \sum_{n=1}^\infty   \frac 1{n^{2\beta-6a}}<\infty.
$$ 
It follows from Proposition~3.4 that $T_*$ is conservative.
It is weakly mixing and of type $III_1$ by Corollary~3.7.
\endexample

We now show that the restriction $a<\frac 1{6}$ in Example~3.8 can not be dropped.
For that we'll need a criterion for dissipativity from \cite{DaKoRo}.

\proclaim{Lemma 3.9}
If $T\in\text{\rom{Aut}}_2(X,\mu)$ and \ 
$
\sum_{n\ge 0}e^{-\frac12\left\|\sqrt{\frac{d\mu\circ T^n}{d\mu}}-1\right\|_2^2}<\infty
$
then $T_*$ is totally dissipative.
\endproclaim

\proclaim{Proposition 3.10} Let $(X,\mu,T)$ be as in Example~3.8 but $a>4$.
Then $T_*$ is totally dissipative.
\endproclaim

\demo{Proof}
For each $n>0$,
$$
\align
\Bigg\|\sqrt{\frac{d\mu\circ T^n}{d\mu}}-1\Bigg\|_2^2&=
a\sum_{k>-n} \left(e^{\frac{\epsilon_{n+k}-\epsilon_k}2}-1\right)^2e^{\epsilon_k}\\
&=
a \sum_{k=-n+1}^{1} \left(e^{\frac{\epsilon_{n+k}}2}-1\right)^2
+a\sum_{k=2}^\infty \left(e^{\frac{\epsilon_{n+k}-\epsilon_k}2}-1\right)^2e^{\epsilon_k}
\\
&=a\sum_{k=1}^{n+1} \left(e^{\frac{\epsilon_{k}}2}-1\right)^2
+a\sum_{k=2}^\infty \left(e^{\frac{\epsilon_{n+k}}2}-e^{\frac{\epsilon_k}{2}}\right)^2\\
&= \frac a4\sum_{k=1}^{n+1} \left(\frac{1}{\sqrt k}+O\bigg(\frac 1k\bigg)\right)^2
+\frac a4\sum_{k=2}^{+\infty}\bigg(\frac1{\sqrt{n+k}}-\frac1{\sqrt{ k}}+O\bigg(\frac1k\bigg)\bigg)^2\\
&= \frac a4\sum_{k=1}^{n+1}\frac1k+\frac a4\sum_{k=2}^{+\infty}\bigg(\frac1{n+k}+\frac1k-\frac2{\sqrt{k(n+k)}}\bigg)+O(1)
\endalign
$$
as $n\to\infty$.
Since $\sum_{k=1}^{n+1}\frac1k=\log n+O(1)$ and
$$
\align
\sum_{k=2}^{+\infty}\bigg(\frac1{n+k}+\frac1k&-\frac2{\sqrt{k(n+k)}}\bigg)=
\int_{2}^{+\infty}\bigg(\frac 1{n+t}+\frac 1t-\frac 2{\sqrt{t(n+t)}}\bigg)\,dt+o(1)\\
&=\bigg(\log(t(n+t))-2\log\bigg(t+\frac n2+\sqrt{t(n+t)}\bigg)\bigg)\bigg|_2^{+\infty}+o(1)\\
&=\log n+ O(1),
\endalign
$$
we obtain that  $\Big\|\sqrt{\frac{d\mu\circ T^n}{d\mu}}-1\Big\|_2^2=\frac a2\log n+O(1)$ as $n\to\infty$.
It follows that  there is a real $D$   such that for each $n>0$,
$$
\left\|\sqrt{\frac{d\mu\circ T^n}{d\mu}}-1\right\|_2^2
\ge D+\frac{a\log n}2.\tag3-14
$$
Since $a>4$,  we deduce from~\thetag{3-14} that
$$
\sum_{n\ge 0}e^{-\frac12\left\|\sqrt{\frac{d\mu\circ T^n}{d\mu}}-1\right\|_2^2}\le\sum_{n\ge 0}
e^{-\frac12(D+\frac{a\log n}2)}=
e^{-\frac D2}
\sum_{n\ge 0}
n^{-\frac a4}<\infty.
$$
It follows now from Lemma~3.9 that $T_*$ is totally dissipative.
\qed
\enddemo

\head 4. Poisson suspensions of type $III_1$ over a non-atomic  base
 \endhead
In this subsection we construct concrete examples of weakly mixing Poisson suspensions of type $III_1$ which are Poisson suspensions of nonsingular transformations defined on non-atomic $\sigma$-finite spaces.

Let $(X,\mu,T)$ be as in \S\,3, i.e. $X=\Bbb Z$, $\mu$  a non-degenerated measure on $X$
and $T$ the unit translation on $X$.
Let $K=[0,1]$. 
Denote by $\lambda$ the Lebesgue measure on $K$.
Let
 $Y:=K\times \Bbb Z$ and let  $Q:Y\to Y$ be the direct product  
 of the identity on $K$ and  $T$
on $X$, i.e.
$Q(k,n):=(k,n+1)$ for all $k\in K$ and $n\in\Bbb Z$.
Then the product measure $\nu:=\lambda\otimes\mu$ is a non-atomic measure
 on $Y$.
Of course, $Q$ is $\nu$-nonsingular.
Moreover,  
$$
Q'(k,n)=T'(n)=\frac{a_{n-1}}{a_n},
$$
where
$a_n:=\mu(\{n\})$.
The dynamical system $(X,\mu,T)$ is a factor of $(Y,\nu,Q)$.
Moreover, the corresponding projection
$$
\vartheta:Y\ni(k,n)\mapsto n\in X
\tag4-1
$$
 is relatively finite measure preserving (r.f.m.p.).
It is straightforward to verify that $Q\in\text{Aut}_2(Y,\nu)$ if and only if  $T\in\text{Aut}_2(X,\mu)$, i.e. if and only if \thetag{3-1} is satisfied.
In a similar way, $Q\in\text{Aut}_1(Y,\nu)$ if and only if
 $T\in\text{Aut}_1(X,\mu)$.
 Hence if $Q\in\text{Aut}_1(Y,\nu)$ then the limits 
 $\lim_{n\to+\infty}a_n$ and $\lim_{n\to-\infty}a_n$ exist and $\chi(Q)=\chi(T)=\lim_{n\to+\infty}a_n-\lim_{n\to-\infty}a_n$.
 
 To describe the Poisson suspension of $(Y,\nu,Q)$ we first write $Y$ as the ``disjoint'' union
 $Y=\bigsqcup_{n\in\Bbb Z}K\times \{n\}$.
Under the natural identification, we may assume that
$$
Y^*=\prod_{n\in\Bbb Z}(K\times\{n\})^*\quad\text{and}\quad \nu^*=\bigotimes_{n\in\Bbb Z}\nu_n^*,
$$
where $\nu_n:=(\nu\restriction (K\times \{n\}))$.
If we identify naturally $K\times \{n\}$ with $K$  then
$\nu_n$ corresponds to $a_n\cdot\lambda$ for each $n\in\Bbb Z$. 
Hence the Poisson suspension $Q_*$ of $Q$ corresponds to the 2-sided Bernoulli shift on the infinite product space
 $$
 (Y^*,\nu^*)=\left((K^*)^{\Bbb Z},\bigotimes_{n\in\Bbb Z}(a_n\cdot\lambda)^*\right).
 $$
We also recall  that we view $(X^*,\mu^*)$ as the infinite product $ ((\Bbb Z_+)^\Bbb Z,\bigotimes_{n\in\Bbb Z}\kappa_n)$, where
$\kappa_n:=(\mu\restriction\{n\})^*$  (see \S\, 3).
The factor map $\vartheta$ from \thetag{4-1} generates a factor map $\vartheta^*:Y^*\to X^*$.
Thus  we have that
 $\vartheta^* Q_*=T_*\vartheta^*$ and $\mu^*=\nu^*\circ(\vartheta^*)^{-1}$.
Since $\vartheta$ is r.f.m.p., $\vartheta^*$ is r.f.m.p. whenever $Q_*$ is $\nu^*$-nonsingular.
It is important to note that  $\vartheta^*$ has a ``coordinate-wise'' structure, i.e.
 $\vartheta^*y=(\vartheta^*_ny_n)_{n\in\Bbb Z}$ for each $y=(y_n)_{n\in\Bbb Z}\in Y^*$, where
$\vartheta^*_n$ is a measure preserving mapping from the probability space $(K^*,\nu_n^*)$ onto the probability space $(\Bbb Z_+,\kappa_n)$ given by $K^*\ni y_n\mapsto y_n(K)\in\Bbb Z_+$.

\proclaim{Corollary 4.1} The Bernoulli shift $Q_*$  is $\nu^*$-nonsingular if and only if $T_*$ is $\mu^*$-nonsingular, i.e. if and only if  \thetag{3-1} is satisfied. 
In this case, the mapping
$$
\vartheta^*:Y^*\ni y=(y_n)_{n\in\Bbb Z}\mapsto (\vartheta^*_ny_n)_{n\in\Bbb Z}:=(y_n(K))_{n\in\Bbb Z}\in X^*
$$
intertwines $Q_*$ with $T_*$, maps $\nu^*$ onto $\mu^*$ and 
$$
\frac{d\nu^*\circ Q_*}{d\nu^*}(y)=\frac{d\mu^*\circ T_*}{d\mu^* }(\vartheta^*y)=\prod_{n\in\Bbb Z}\frac{\kappa_{n-1}(\vartheta_n^*y_n)}{\kappa_n(\vartheta_n^*y_n)}
$$
for a.e. $y=(y_n)_{n\in\Bbb Z}\in Y^*$.
Hence $(Y^*,\nu^*)$ is isomorphic to the product space $(K^\Bbb Z\times X^*,\lambda^\Bbb Z\times\mu^*)$ in such a way that 
\roster
\item"\rom{(i)}" $\vartheta^*$ corresponds to the projection to the second coordinate and 
\item"\rom{(ii)}" $Q_*$ corresponds the direct product $B\times T_*$, where $B:K^\Bbb Z\to K^\Bbb Z$ is the Bernoulli shift preserving $\lambda^\Bbb Z$.
\endroster
\endproclaim

\demo{Proof}
For each $n\in\Bbb Z$, we disintegrate $\nu_n^*$ relative to $\kappa_n$ (via $\vartheta^*_n$):
$$
\nu_n^*=\sum_{k\in\Bbb Z_+}\kappa_n(k)\xi_{n,k},\tag4-2
$$
where $(\xi_{n,k})_{k\in\Bbb Z_+}$ is the corresponding sequence of conditional probabilities on $K^*$.
Using \thetag{4-2}, we obtain the disintegration of $\nu^*$ relative to $\mu^*$ (via $\vartheta^*$):
$$
\nu^*=\int_{X^*}\bigotimes_{n\in\Bbb Z}\xi_{n,k_n}\,d\mu^*((k_n)_{n\in\Bbb Z}).
$$
Since $\vartheta^*$ intertwines the Bernoulli shifts on $Y^*$ and $X^*$ respectively, the r.f.m.p. property of $\vartheta^*$ means that
$
\bigotimes_{n\in\Bbb Z}\xi_{n,k_{n+1}}=\bigotimes_{n\in\Bbb Z}\xi_{n+1,k_{n+1}}
$
and hence $\xi_{n,k_{n+1}}=\xi_{n+1,k_{n+1}}$
for all $n\in\Bbb Z$ and $\mu^*$-a.e. $(k_n)_{n\in\Bbb Z}\in X^*$.
Hence there is a sequence $(\xi_k)_{k\in\Bbb Z_+}$ of probability measures on $K^*$ such that 
$\xi_{n,k}=\xi_k$ 
and therefore
$$
\nu^*=\int_{X^*}\bigotimes_{n\in\Bbb Z}\xi_{k_n}\,d\mu^*((k_n)_{n\in\Bbb Z}).
$$ 
The measure $\xi_k$ is non-atomic for each $k>0$ and the measure $\xi_0$ is a delta-measure.
For a.e. $\boldkey k=(k_n)_{n\in\Bbb Z}\in X^*$, there exist infinitely many positive $n$ with $k_n>0$ and infinitely many negative $n$ with $k_n>0$.
Hence for such a $\boldkey k$ and every  $n\in\Bbb Z$, there are (uniquely defined) two integers $l_n(\boldkey k)$ and $r_n(\boldkey k)$ such that the following are satisfied:
\roster
\item"---"
$l_n(\boldkey k)\le n\le r_n(\boldkey k)$,
\item"---"
$k_{l_n(\boldkey k)}\ne 0$ and $k_{r_n(\boldkey k)+1}\ne 0$ and
\item"---"
 if $l_n(\boldkey k)< n\le r_n(\boldkey k)$ then  $k_n=0$.
\endroster
Of course, the mappings  $X^*\ni \boldkey k\mapsto l_n(\boldkey k)\in\Bbb Z$ and 
$X^*\ni \boldkey k\mapsto r_n(\boldkey k)\in\Bbb Z$ are measurable for every $n\in\Bbb Z$.
For each $n\in\Bbb Z\setminus\{0\}$ and $l\in\Bbb N$,
there is  a measure theoretical isomorphism 
$$
\tau_{k,l}:((K^*)^l,\xi_k\otimes\xi_0\otimes\cdots\otimes\xi_0)\to(K^l,\lambda^{ l}).
$$
Next, for each point $\boldkey k=(k_n)_{n\in\Bbb Z}\in X^*$ and $\boldkey y=(y_n)_{n\in\Bbb Z}\in (K^*)^\Bbb Z$, we define
a point $\langle{\boldkey k},\boldkey y\rangle=(z_n)_{n\in\Bbb Z}$ of $K^\Bbb Z$ in the following way:
$z_n$ is the $(n-l_n(\boldkey k)+1)$-th symbol in the block  
$$
\tau_{l_n(\boldkey k), r_n(\boldkey k)-l_n(\boldkey k)+1}[y_{l_n}\cdots y_{r_n}]\in K^{r_n(\boldkey k)-l_n(\boldkey k)+1}.
$$
We now define an isomorphism $\tau$ of $(Y^*,\nu^*)$ onto $(K^\Bbb Z\times X^*,\lambda^\Bbb Z\otimes\mu^*)$ by setting
$$
\tau(\boldkey y):=(\langle  \vartheta^*\boldkey y,\boldkey y  \rangle,
\vartheta^*\boldkey y).
$$
It is straightforward to verify that (i) and (ii) hold. \qed

\enddemo

We now state the main result of this section.

\proclaim{Theorem 4.2} If \thetag{3-1} and \thetag{3-4} are satisfied then
 the Maharam extension of $(Y^*,\nu^*,Q_*)$  is a  $K$-automorphism.
 If, moreover, $T_*$ is conservative then $Q_*$ is weakly mixing and of stable Krieger's type $III_1$.
\endproclaim

\demo{Proof}
Consider $Q_*$ as the direct product $B\times T_*$ described in Corollary~4.1.
Then the Maharam extension of $Q_*$ is isomorphic to the product of $B$ with the Maharam extension of  $T_*$.
It remains to apply Theorem~3.5 and a simple fact that the direct product of a (conservative) $K$-automorphism with a probability preserving Bernoulli shift is a (conservative) $K$-automorphism.
\qed

\enddemo

\remark{Remark \rom{4.3}} Utilizing  Example 3.8  we obtain concrete examples of ergodic conservative suspensions $Q_*$ of type $III_1$  for $Q$ defined on a nonatomic infinite measure space.
Moreover, it follows from Example~3.8 and Proposition~3.10 that $Q_*$ is totally dissipative if $a>4$ and conservative (and hence weakly mixing and of stable type $III_1$) if $0<a<\frac1{6}$.
This result will be improved in \S\,7 below.
\endremark

\head 5. Generic Poisson suspension is of type $III_1$
 \endhead
Throughout this section we assume that $\mu$ is an infinite $\sigma$-finite non-atomic measure
 on standard Borel space $(X,\goth B)$.
Given two subsets $A,B\subset X$ of finite measure, let $\tau_{A,B}$ denote a $\mu$-nonsingular bijection from $A$ onto $B$ such that  $\frac{d\mu\circ\tau_{A,B}}{d\mu}(x)=\frac{\mu(B)}{\mu(A)}$ for all $x\in A$.
 We now recall a definition and some facts from \cite{DaKoRo}.
 
 \definition{Definition 5.1}
A nonsingular transformation $T$ of a $\sigma$-finite standard measure space $(X,\mu)$ is {\it locally aperiodic} if there is a subset $A\subset X$ of positive finite measure such that $Tx=x$ if $x\not\in A$ and
 $T^nx\ne x$ for each $x\in A$ and $n>0$.
 \enddefinition
 
 Of course, each locally aperiodic transformation $T$ belongs to Aut$_1(X,\mu)$.
It is easy to verify that $\chi(T)=0$. 
 It was also shown in \cite{DaKoRo} that:
\roster
\item"(F1)"
the conjugacy class of $T$ in  $\text{\rom{Aut}}_2(X,\mu)$ is $d_2$-dense in $\text{\rom{Aut}}_2(X,\mu)$ and
\item"(F2)"
the conjugacy class of $T$ in  $\text{\rom{Aut}}_1(X,\mu)$ is $d_1$-dense in $\text{\rom{Ker\,}}\chi$.
\endroster

 \proclaim{Proposition 5.2} Let  $T\in\text{\rom{Aut}}(X,\mu)$ and there are a partition
$X=\bigsqcup_{n\in\Bbb Z}W_n$ of $X$ into subsets $W_n$ of finite measure and a sequence $(a_n)_{n\in\Bbb Z}$ of reals such that $TW_n=W_{n+1}$ and $T'(x)=a_n$ at a.e. $x\in W_n$ and  for each $n\in\Bbb Z$.
Suppose that there exists a limit $\lim_{|n|\to\infty}\mu(W_n)\in(0,+\infty)$.
\roster
\item"\rom(i)" If $T\in\text{\rom{Aut}}_2(X,\mu)$ then
 the conjugacy class of $T$ is dense in $\text{\rom{Aut}}_2(X,\mu)$.
 \item"\rom{(ii)}"
If $T\in\text{\rom{Aut}}_1(X,\mu)$ then
 the conjugacy class of $T$ is dense in $\text{\rom{Ker}}\,\chi$.
\endroster
\endproclaim

\demo{Proof} (i)
Since $T\in\text{\rom{Aut}}_2(X,\mu)$, 
$$
\infty>\|\sqrt{T'}-1\|_2^2=\sum_{n\in\Bbb Z}\int_{W_n}(\sqrt{T'(x)}-1)^2d\mu(x)=
\sum_{n\in\Bbb Z}(\sqrt{a_n}-1)^2\mu(W_n).
$$ 
Since 
 $a_n=\mu(W_{n+1})/\mu(W_n)$, we obtain that
$$
\sum_{n\in\Bbb Z}(\sqrt{\mu(W_{n+1})}-\sqrt{\mu(W_n)})^2<\infty.\tag5-1
$$
If we show that the $d_2$-closure of the conjugacy class of $T$ contains a locally aperiodic transformation then (F1) yields that the conjugacy class of $T$ is dense in Aut$_2(X,\mu)$.
We now construct such a transformation.
For each $n>0$, we set  $h_n:=2^n$.
Let $B_0$ be a subset of $X$ with $\mu(B_0)=1$ and let $S$ be a $\mu$-preserving transformation of $X$ such that $Sx=x$ if $x\not\in B_0$ 
and the restriction of $S$ to $B_0$ is isomorphic to the  2-adic odometer.
Then $S$ is locally aperiodic and
 there exists a decreasing infinite sequence of subsets $B_0\supset B_1\supset\cdots$ in $X$
such that
$S^iB_n\cap S^jB_n=\emptyset$ if $0\le i<j<h_n$ and $\bigsqcup_{i=0}^{h_n-1}S^iB_n=B_0$.
We see that $S$ has a ``cyclic structure of period $h_n$ on $B_n$".
The idea  of the following argument is to define for each $n>0$, a transformation $S_n$
extending   the cyclic structure of $S$ from $B_n$ to a larger subset
in such a way $S_n\to S$ and $S_n$ is conjugate to a cyclic permutation of the finite sequence $W_{-h_{n-1}},\dots, W_{h_{n-1}-1}$ that approaches to $T$ as $n\to\infty$:
Since the convergence is considered in $d_2$, we should  choose carefully the transformations in these constructions to control their Radon-Nikodym derivatives.

More precisely, we now construct a sequence of $\mu$-preserving transformations $(S_n)_{n=1}^\infty$ 
and a sequence $(B_n')_{n=1}^\infty$ of subsets in $X$ such that
\roster
\item"---"
$S_n\restriction B_0=S$ for each $n$ and
$S_n\to S$ weakly (and hence in $d_2$) as $n\to\infty$, 
\item"---"$
B_n'\supset B_n$, $\mu(B_n')=\mu(W_{-h_{n-1}})$, 
$S^i_nB_n'\cap S^j_nB_n'=\emptyset$ for all $0\le i<j<h_n$ and each $n$.
\item"---"
If we set $Y_n:=\bigsqcup_{i=0}^{h_n-1}S_n^iB_n'$
and $Y_n^\circ:=Y_n\setminus B_n'$ then for each subset $A\subset X$ of finite measure,
$\lim_{n\to\infty}\mu(A\cap Y_n^\circ)=\mu(A)$.
\endroster
To see that such a construction is possible, it is convenient to think that $X=[0,+\infty)$, $\mu$ is a Lebesgue measure on $X$, $B_0=[0,1)$, $B_n:=[0,h_n^{-1})$.
To define $B_n'$ and $S_n$, we first select  a real $\delta_n>0$ and an integer $D_n>0$ such that
\roster
\item"(a)" $\delta_n h_n<n^{-1}$ and
\item"(b)" $\mu(B_n)+\delta_n D_n=\mu(W_{-h_{n-1}})$.
\endroster
Now we 
set $B_n':=B_n\sqcup\bigsqcup_{j=0}^{D_n-1}[1+jh_n\delta_n,1+(jh_n+1)\delta_n)$ and
$$
S_nx:=
\cases
Sx, &\text{if } x\in B_0\\
x+\delta_n, &\text{if }1\le x <1+(D_nh_n-1)\delta_n\\
x-(D_nh_n-1)\delta_n , &\text{if }1+(D_nh_n-1)\delta_n\le x<1+D_nh_n\delta_n\\
x, &\text{if }x\ge 1+D_nh_n\delta_n.
\endcases
$$
It follows from (a) that $S_n\to S$ weakly as $n\to\infty$.
(b) implies that $\mu(B_n')=\mu(W_{-h_{n-1}})$.
The sets $S_n^iB_n'$, $i=0,\dots,h_n-1$, are mutually disjoint and their union $Y_n$ equals $[0,1+h_nD_n\delta_n)$.
In particular, $Y_n\supset [0,1+h_nb/2)$ eventually in $n$, where $b:=\lim_{|m|\to\infty}\mu(W_m)$.
Hence $\bigcup_{n\to\infty}Y_n=X$.
Since the set $B_n'\setminus B_n$ is ``uniformly distributed'' along $Y_n\setminus B_0$, we have that $\mu(Y_n^\circ\cap [0,E))\to\mu([0,E))$ for each $E>1$.
Thus, the sequences $(S_n)_{n=1}^\infty$ and $(B_n')_{n=1}^\infty$
are as desired.

We now set
$X_n:=\bigsqcup_{i=0}^{h_{n}-1}T^iW_{-h_{n-1}}$
and $X_n^\circ:=X_n\setminus W_{-h_{n-1}}$.
Then
 $X_1^\circ\subset X_2^\circ\subset\cdots$ and $\bigcup_{n=0}^\infty X_n^\circ=X$.
For each $n>0$, we select a measure preserving Borel bijection $\tau_n$ of $X\setminus X_n$
onto $X\setminus Y_n$
and define
a transformation $R_n$ of $X$ by setting
$$
R_nx:=
\cases
S_n^i\tau_{W_{-h_{n-1}},B_n'}T^{-i}&\text{ if $x\in W_{i-h_{n-1}}$, $0\le i<h_n$.}\\
\tau_nx&\text{if $x\not\in X_n$}.
\endcases
$$
Then $R_n\in\text{Aut}(X,\mu)$, $R_nX_n=Y_n$ and 
hence $R_n\in\text{Aut}_1(X,\mu)$.
A straightforward verification shows that
$R_nT^{-1}R_n^{-1}=S_n^{-1}$ on $Y_n^\circ$.
Therefore $S_nR_nT^{-1}R_n^{-1}\to\text{Id}$ weakly as $n\to\infty$.
Next, we note that the Radon-Nikodym derivatives of  $T^{-1}$ and $R_n$ are constant on each level $W_j$, $j\in\Bbb Z$.
Since $S_n$ preserves $\mu$ and $T^{-1}W_j=W_{j-1}$ for each $j\in\Bbb Z$, it follows that    the Radon-Nikodym derivative of the transformation
$S_nR_nT^{-1}R_n^{-1}$ is constant on the subset  $R_nW_j$ for each $j\in\Bbb Z$.
Hence we compute easily that $(S_nR_nT^{-1}R_n^{-1})' (x)$ equals
$$
\cases
1,&\text{if }x\in Y_n^\circ \\
\frac{\mu(S_nR_nW_{i-h_{n-1}-1})}{\mu(R_nW_{i-h_{n-1}})},
& \text{if }x\in  R_nW_{i-h_{n-1}}\text{ and }i\not\in\{1,\dots,h_n-1\}.
\endcases
$$
If $i\not\in\{1,\dots,h_n\}$ then $\mu(S_nR_nW_{i-h_{n-1}-1})=\mu(W_{i-h_{n-1}-1})$.
If $i=h_n$ then
$\mu(S_nR_nW_{i-h_{n-1}-1})=\mu(S_nS_n^{h_{n}-1}\tau_{W_{-h_{n-1}},B_n'}W_{-h_{n-1}})=\mu(W_{-h_{n-1}})$.
Therefore
$$
(S_nR_nT^{-1}R_n^{-1})' (x)  
=
\cases
1,&\text{if }x\in Y_n^\circ\\
\frac{\mu(W_{i-h_{n-1}-1})}{\mu(W_{i-h_{n-1}})},&
 \text{if }x\in  R_nW_{i-h_{n-1}}\text{ and }i\not\in\{1,\dots,h_n\}\\
\frac{\mu(W_{-h_{n-1}})}{\mu(W_{h_n-h_{n-1}})},& \text{if }x\in R_nW_{h_n-h_{n-1}}.
\endcases
$$
Therefore
$$
\aligned
\left\|\sqrt{(S_nR_nT^{-1}R_n^{-1})'}-1\right\|_2^2&=\sum_{i\not\in\{1,\dots,h_n\}}\left(\sqrt{\mu(W_{i-h_{n-1}-1})}-\sqrt{\mu(W_{i-h_{n-1}})}\right)^2\\
&+
\left(\sqrt{\mu(W_{-h_{n-1}})}-\sqrt{\mu(W_{h_n-h_{n-1}})}\right)^2.
\endaligned
$$
Since  $b<\infty$, it follows that $\sqrt{\mu(W_{-h_{n-1}})}-\sqrt{\mu(W_{h_{n-1}})}\to 0$ as $n\to\infty$.
Utilizing  this fact and \thetag{5-1} we obtain that
$$
\lim_{n\to\infty}\left\|\sqrt{(S_nR_nT^{-1}R_n^{-1})'}-1\right\|_2= 0.
$$
Therefore $S_nR_nT^{-1}R_n^{-1}\to \text{Id}$ in $d_2$ as $n\to\infty$.
Since $S_n\to S$ in $d_2$, we obtain that
 $SR_nT^{-1}R_n^{-1}\to \text{Id}$ in $d_2$ as $n\to\infty$.
Thus $S$ belongs to the $d_2$-closure of the conjugacy class of $T$ in  Aut$_2(X,\mu)$, as desired.

(ii) is proved in a similar way. 
First, we note that the conditions $T\in\text{Aut}_1(X,\mu)$ and  $\lim_{|n|\to\infty}\mu(W_n)<\infty$ imply that $T\in\text{Ker}\,\chi$.
Indeed,
$$
\chi(T)=\int_X(T'-1)\,d\mu=\sum_{n\in\Bbb Z}\int_{W_n}(T'-1)\,d\mu=\sum_{n\in\Bbb Z}(\mu(W_{n-1})-\mu(W_n))=0.
$$
Secondly, we will use (F2) in place of (F1).
Instead of \thetag{5-1} we now have that
$\sum_{n\in\Bbb Z}|{\mu(W_{n+1})}-{\mu(W_n)}|<\infty$.
Define $R_n$ in the same way as in the proof of (i).
As was noted there, $R_n\in\text{Aut}_1(X,\mu)$.
Moreover, 
it is easy to see that  $R_n\in \text{Ker}\,\chi$.
Slightly modifying the above argument and  considering the $L^1$-norm  instead of  the square of the $L^2$-norm, we obtain that $\|{(S_nR_nT^{-1}R_n^{-1})'}-1\|_1\to 0$.
The latter yields that $SR_nT^{-1}R_n^{-1}\to \text{Id}$ in $d_1$ as $n\to\infty$.
\qed
\enddemo

We recall two more facts from \cite{DaKoRo}.
Let $\Cal E$ denote the subset of all ergodic transformations of type $III_1$ in Aut$(X,\mu)$. 
\roster
\item"(F3)"  $\Cal E\cap\text{Aut}_2(X,\mu)$ is a dense $G_\delta$ in $(\text{Aut}_2(X,\mu),d_2)$.
\item"(F4)" $\Cal E\cap\text{Aut}_1(X,\mu)=\Cal E\cap\text{Ker\,}\chi$ is a dense $G_\delta$ in $(\text{Ker\,}\chi,d_1)$.
\endroster

We now state and prove the main result of the section.

\proclaim{Theorem  5.3} 
\roster
\item"\rom{(i)}" The subset
$
\Cal E_2^*:=\{T\in\Cal E\cap \text{\rom{Aut}}_2(X,\mu)\mid \text{ $T_*$  is ergodic and of type $III_1$}\}
$
 is a dense $G_\delta$ in {\rom{(Aut}}$_2(X,\mu),d_2)$. 
 \item"\rom{(ii)}"
 The subset
$
\Cal E_1^*:=\{T\in\Cal E\cap \text{\rom{Ker}}\,\chi\mid \text{ $T_*$  is ergodic and of type $III_1$}\}
$
 is a dense $G_\delta$ in  $\text{\rom{(Ker}}\,\chi,d_1)$. 

 \endroster
\endproclaim
\demo{Proof}
(i) Since the set of ergodic transformations of type $III_1$  is $G_\delta$ in Aut$(X^*,\mu^*)$ endowed with the weak topology \cite{ChHaPr} and the
 map 
 $$
 \text{Aut}_2(X,\mu)\ni T\mapsto T_*\in\text{Aut}(X^*,\mu^*)
 $$
  is continuous \cite{DaKoRo}, it follows   that the subset
  $$
  \Cal E_2:=\{T\in\text{Aut}_2(X,\mu)\mid T_*\text{ is ergodic of type $III_1$}\}
  $$
is a $G_\delta$ in (Aut$_2(X,\mu),d_2)$. 
 By Theorem~4.2 and Example~3.8,  $\Cal E_2\ne\emptyset$.
Moreover, 
 there exists a transformation $T\in \Cal E_2$ satisfying the conditions of Proposition~5.2.
 Hence, by Proposition~5.2, the conjugacy class of $T$ is dense in (Aut$_2(X,\mu),d_2)$. 
 Since $\Cal E_2$ is conjugacy invariant, it follows that $\Cal E_2$ is a dense $G_\delta$.
Of course, $\Cal E_2^*=\Cal E_2\cap \Cal E\cap\text{Aut}_2(X,\mu)$.
We  finally deduce  (i) from these facts and  \thetag{F3}.

(ii) is proved in a similar way with usage of (F4) instead of (F3). 
 \qed
\enddemo

\head 6. Poisson suspensions over  dissipative bases
\endhead
Let $(X,\goth B,\mu)$ be a $\sigma$-finite infinite standard nonatomic  measure space and let
$T\in\text{Aut}_2(X,\mu)$.
Let $T$ be {\it totally dissipative}.
This means that there is a subset $B\in\goth B$ such that $X=\bigsqcup_{n\in\Bbb Z}T^nB$.
Suppose  that $\mu(B)<\infty$.
Since $T\in\text{Aut}_2(X,\mu)$, it follows $\mu(T^nB)<\infty$ for each $n\in\Bbb Z$ (see \cite{DaKoRo}).
Then without loss of generality we can assume that $X=[0,1]\times\Bbb Z$ and $T(y,n)=(y,n+1)$ for all $y\in[0,1]$ and  $n\in\Bbb Z$.
Let $\lambda$ denote the Lebesgue measure on $[0,1]$.
We can also assume that
there exist a sequence $(a_n)_{n\in\Bbb Z}$ of  functions $a_n\in L^1([0,1],\lambda)$ such that $a_n>0$ a.e. and for each $F\in L^1([0,1]\times\Bbb Z,\mu)$,
$$
\int_X Fd\mu=\sum_{n\in\Bbb Z}\int_{[0,1]}F(y,n)a_n(y)\,d\lambda(y).\tag6-1
$$ 
It is straightforward to verify that  $T\in\text{Aut}_2(X,\mu)$ if and only if
$$
\sum_{n\in\Bbb Z}\int_{[0,1]}\left(\sqrt{\frac{a_{n+1}(y)}{a_n(y)}}-1\right)^2a_n(y)\,d\lambda(y)=
\sum_{n\in\Bbb Z}\big\|\sqrt{a_{n+1}}-\sqrt{a_n}\big\|_2^2<\infty.
$$
In a similar way, $T\in\text{Aut}_1(X,\mu)$ if and only if $\sum_{n\in\Bbb Z}\|{a_{n+1}}-{a_n}\|_1<\infty.$
The latter inequality  implies that there exist two nonnegative functions $a,b\in L^1([0,1],\lambda)$
such that 
$$
\lim_{n\to\infty}\|a-a_n\|_1=0\quad\text{ and }
\lim_{n\to-\infty}\|b-a_n\|_1=0.\tag6-2
$$
Since $\chi(T)=\lim_{m\to\infty}(\mu(T([0,1]\times\{-m,\dots,m\}))-\mu([0,1]\times\{-m,\dots,m\}))$, it follows that
$$
\chi(T)=\lim_{m\to\infty}\left(\sum_{j=-m+1}^{m+1}\int a_jd\lambda-\sum_{j=-m}^m\int a_jd\lambda\right)=
\|a\|_1-\|b\|_1.
\tag6-3
$$
\comment

We will need the following lemma.

\proclaim{Lemma 6.1} Let $\kappa$ and $\nu$ be two finite measures on $[0,1]$ such that
$\kappa\sim\nu\sim\lambda$.
Then
$$
-2\log\int_{[0,1]^*}\sqrt{\frac{d\kappa^*}{d\lambda^*}\frac{d\nu^*}{d\lambda^*}}d\lambda^*=
\left\|\sqrt{\frac{d\kappa}{d\lambda}}-\sqrt{\frac{d\nu}{d\lambda}}\right\|_2^2.
$$
\endproclaim

\endcomment

The following theorem is the main result of this section.

\proclaim{Theorem 6.1} If  $T\in\text{\rom{Aut}}_1(X,\mu)$ and $\chi(T)\ne 0$ then $T_*$ is totally dissipative.
\endproclaim
\demo{Proof} 
Consider the Hopf decomposition of $X$:  there are two invariant subsets  $X_c$ and $X_d$
of $X$ such that the dynamical system $(X_c,\mu\restriction X_c, T\restriction X_c)$ is conservative and the dynamical system $(X_d,\mu\restriction X_d, T\restriction X_d)$ is  totally dissipative \cite{Aa}.
Since $\chi(T)\ne 0$, it follows that $T$ is dissipative \cite{DaKoRo} and 
hence $\mu(X_d)>0$.
Of course, $T\restriction X_d\in\text{\rom{Aut}}_1(X_d,\mu_d)$, $T\restriction X_c\in\text{\rom{Aut}}_1(X_c,\mu_c)$ and 
$$
\chi(T)=\chi(T\restriction X_d)+\chi(T\restriction X_c)=\chi(T\restriction X_d).
$$
Hence $\chi(T\restriction X_d)\ne 0$.
Moreover, $T_*$ is isomorphic to the Cartesian product of $(T\restriction X_d)_*$ and
$(T\restriction X_c)_*$.
Therefore if $(T\restriction X_d)_*$ is totally dissipative then so is $T_*$.
In view of this, we may assume without loss of generality that $T$ is totally dissipative itself.
Hence there is a subset $B\subset X$ such that
$X=\bigsqcup_{n\in\Bbb Z}T^nB$.
We may also assume that $\mu(B)<\infty$.
Indeed, if $\mu(B)=\infty$ then we partition $B$ into countably many subsets $B_m$, $m\in\Bbb N$ of finite measure.
For each $m>0$, let $B_m':=\bigsqcup_{n\in\Bbb Z}T^nB_m$.
Then $B_m'$ is invariant under $T$ and $\bigsqcup_{m=1}^\infty B_m'=X$.
Of course, $T\restriction B_m'\in\text{Aut}_1(B_m',\mu\restriction B_m')$ for each $m>0$ and
$\chi(T)=\sum_{m=1}^\infty \chi(T\restriction B_m')$.
Therefore there exists $m_0>0$ such that $\chi(T\restriction B_{m_0}')\ne 0$.
Since $T_*$ is isomorphic to the Cartesian  product of $(T\restriction B_{m_0}')_*$ and
$(T\restriction (X\setminus B_{m_0}'))_*$, it follows that if $(T\restriction B_{m_0}')_*$ is totally dissipative then so is $T_*$.
Thus, it suffices to consider only the case where $\mu(B)<\infty$.

In this case we can consider $X$ as $[0,1]\times\Bbb Z$ and $T$ as the unit rotation along the second coordinate.
 We will also use the notation $(a_n)_{n\in\Bbb Z}$  and $\lambda$ introduced in the beginning of this section  to describe $\mu$ via \thetag{6-1}.
 Then $T_*$ is the Bernoulli shift on the infinite product space $X^*:=\bigotimes_{n\in\Bbb Z}([0,1]^*,\kappa_n^*)$, where $\kappa_n$ is a measure equivalent to $\lambda$ on $[0,1]$ with
 $\frac{d\kappa_n}{d\lambda}=a_n$.
 We are going to apply Lemma~3.9 to prove that $T_*$ is totally dissipative.
We will use the following equality:
 $$
 \aligned
 \left\|\sqrt{\frac{d\mu\circ T^n}{d\mu}}-1\right\|_2^2&=\sum_{k\in\Bbb Z}\int_{[0,1]}\left(\sqrt{\frac{a_{k+n}}{a_k}}-1\right)^2a_kd\lambda\\
 &=\sum_{k\in\Bbb Z}\left\|\sqrt{{a_{k+n}}}-\sqrt{a_k}\right\|_2^2.
 \endaligned
 \tag6-4
 $$
Since for two nonnegative functions $e,f\in L^1([0,1],\lambda)$, 
$$
\|e-f\|_1=\langle\sqrt e-\sqrt f,\sqrt e+\sqrt f\rangle\le\| \sqrt e-\sqrt f\|_2\cdot \|\sqrt e+\sqrt f\|_2,
$$
and  $\|\sqrt e+\sqrt f\|_2\le \|\sqrt e\|_2+\|\sqrt f\|_2=\sqrt{ \| e\|_1}+\sqrt{\| f\|_1}$,  we obtain  that
 $$
 \left\|\sqrt{{a_{k+n}}}-\sqrt{a_k}\right\|_2\ge \frac 1{2D}\left\|a_{k+n}-a_k\right\|_1\tag6-5
 $$
 for all $k$ and $n$,
 where $D:=\sup_{k\in\Bbb Z}\|a_k\|_1<\infty$.
 \comment

 By \cite{DaKoRo}, $T_*$ is $\mu^*$-nonsingular and
 $$
(T_*)'(y)=\prod_{n\in\Bbb Z}\frac{\kappa^*_{n-1}(y_n)}{\kappa^*_n(y_n)}\quad\text{for a.a. $y=(y_n)_{n\in\Bbb Z}\in X^*$.}\tag6-4
$$

\endcomment
It follows from \thetag{6-2} and \thetag{6-3} that there is  $N>0$ such that 
$$
\|a_i-a_j\|_1>\frac 13|\chi(T)|\quad\text{whenever $i>N$ and $j<-N$.}\tag6-6
$$ 
Taking $n>3N$ and utilizing  \thetag{6-4}, \thetag{6-5}, \thetag{6-6}  
we obtain that 
$$
 \aligned
 \left\|\sqrt{\frac{d\mu\circ T^n}{d\mu}}-1\right\|_2^2&\ge
 \sum_{k\in\Bbb Z}\frac {\|a_{k+n}-a_k\|_1^2}{4D^2}\\
& \ge
 \sum_{ -\frac n3>k>-\frac{2n}3}\frac {\|a_{k+n}-a_k\|_1^2}{4D^2}\\
 &\ge n \frac {|\chi(T)|^2} {108D^2}.
 \endaligned
$$
Therefore 
$$
\sum_{n>0}e^{-\frac 12 \left\|\sqrt{\frac{d\mu\circ T^n}{d\mu}}-1\right\|_2^2}<\sum_{n>0}e^{-\frac {|\chi(T)|^2} {216D^2}n}<\infty.
$$
\comment

$$
\align
\log\langle U_{T_*}^n1,1\rangle &=\log\left(\prod_{n\in\Bbb Z}\int_{[0,1]^*}\sqrt{\frac{d\kappa^*_{k-n}}{d\kappa^*_k}(y)}\,d\kappa^*_k(y)\right)\\
&=
\sum_{k\in\Bbb Z}\log\left(\int_{[0,1]^*}\sqrt{\frac{d\kappa^*_{k-n}}{d\kappa^*_k}(y)}\,d\kappa^*_k(y)\right)\\
&=
-\frac 12\sum_{k\in\Bbb Z}\|\sqrt{a_{k-n}}-\sqrt{a_k}\|_2^2
\\
&\le
-\frac 12\sum_{\frac n3<k<\frac{2n}3}\|\sqrt{a_{k-n}}-\sqrt{a_k}\|_2^2
\\
&\le
-\frac 18\sum_{\frac n3<k<\frac{2n}3}\|{a_{k-n}}-{a_k}\|_1^2
\\
&\le
-\frac n{72} \cdot |\chi(T)|^2.
\endalign
$$
Hence $\sum_{n=1}^\infty\langle U_{T_*}^n1,1\rangle<\sum_{n=1}^\infty e^{-\frac{|\chi(T)|^2}{72}n}
<\infty$. 
\endcomment
 Hence $T_*$ is totally dissipative by  Lemma~3.9.
\qed

\enddemo

\remark{Remark 6.2} 
 This theorem can be extended to the case where $T\in\text{Aut}_2(X,\mu)$ in the following way.
 Suppose that $T$ is dissipative.
  Select a subset $B\subset X$ of finite positive measure  such that $T^nB\cap T^mB=\emptyset$ if $n,m\in\Bbb Z$ and $n\ne m$.
 Then we represent the restriction of $T$ to $\bigsqcup_{n\in\Bbb Z}T^nB$ in the same way as in the beginning of this section such that \thetag{6-1} holds.
 If there exist $\delta>0$ and $N>0$ such that for all $n,m>N$, we have that
 $\|\sqrt{ a_n}-\sqrt{a_{-m}}\|_2>\delta$ then
 $T_*$ is totally dissipative. 
 This fact is proved in the very same way as Theorem~6.1.
 We leave details to the reader.
\endremark

\head Section 7. Phase transition for conservativeness  when scaling the intensity 
of Poisson suspensions
\endhead

Let $(X,\goth B,\mu)$ be a $\sigma$-finite standard measure space with $\mu$ nonatomic and infinite.
For $t>0$, let $\mu_t$ denote the measure  on $(X,\goth B)$ given by $\mu_t(B):=t\mu(B)$.

\definition{Definition 7.1} Let $T\in\text{\rom{Aut}}_2(X,\mu)$.
We say that $T_*$ is {\it conservatively concrete for intensity scaling  (CCIS)}  if for each $t>0$, 
the Poisson suspension $(X^*,\goth B^*,\mu_{t}^*,T_*)$ is  either conservative or totally dissipative.
\enddefinition

If \, $T\in\text{\rom{Aut}}_1(X,\mu)$ and $\,\chi( T)>0$ then $T$ is CCIS.
On the other hand, if $T\in\text{\rom{Aut}}_2(X,\mu)$
and there is a $T$-invariant subset $A\subset X$ of finite positive measure such that
$T\restriction A$ is totally dissipative and $T\restriction (X\setminus A)$ preserves  $\mu\restriction(X\setminus A)$ then $T_*$ is not CCIS.
Indeed, $(X^*,\mu^*,T_*)$ is isomorphic to the Cartesian product 
$$
(A^*,(\mu\restriction A)^*,(T\restriction A)_*)\times((X\setminus A)^*,(\mu\restriction(X\setminus A))^*,(T\restriction (X\setminus A))_*).
$$
We now set $[A]_j:=\{\omega\in X^*\mid \omega(A)=j\}$ for $j=0,1$.
Then $[A]_j$ is invariant under $T_*$ and $\mu^*([A]_j)>0$ for  $j=0,1$.
Since $[A]_0$ is a singleton (modulo $(\mu\restriction A)^*$), it follows that  $T_*$ is conservative when restricted to this subset.
One can check that  $T_*\restriction [A]_1$ is totally dissipative.
Since $(T\restriction (X\setminus A))_*$ is conservative (because it is probability preserving),
it follows that $(T\restriction A)_*\times (T\restriction (X\setminus A))_*$ is conservative on $[A]_0\times(X\setminus A)^*$ and dissipative
 on $[A]_1\times(X\setminus A)^*$.
Hence $T_*$ is not CCIS.

\remark{Problem} What are necessary and sufficient conditions on $(X,\mu,T)$ under which $T_*$ is CCIS?  
\endremark

In  this paper we prove CCIS for a certain  family of Poisson suspensions.
In the proof of the following theorem we use an idea similar to what was used in \cite{Ko2} (see also \cite{Da}).

\proclaim{Theorem 7.2} Let $T\in\text{\rom{Aut}}_1(X,\mu)$.
Suppose also that there is  $\alpha>1$ such that $\alpha^{-1}<(T^n)'(x)<\alpha$ for each $n>0$ at a.e. $x\in X$.
Then $T_*$ is either conservative or totally dissipative.
Hence $T_*$ is CCIS.
\endproclaim

\demo{Proof} 
In view of the aforementioned remark, it suffices to consider only the case where $\chi(T)=0$.

Let $\Cal S$ be the subgroup of all $\mu$-preserving transformations on $(X,\mu)$ such that for each $S\in \Cal S$, there is a subset $A_S\subset X$ of finite measure with $Sx=x$ whenever $x\not\in A_S$.
Let
$
D:=\big\{\omega\in X^*\,\big|\, \sum_{n=1}^\infty (T_*^n)'(\omega)<\infty\big\}.
$

We first show that $D$ is invariant under $S_*$  for each $S\in\Cal S$.
Indeed, for each $n>0$,  $\omega\in D$ and $S\in\Cal S$, we have (utilizing  \cite{DaKoRo})
that
$$
\align
(T_*^n)'(S_*\omega)&=\prod_{\omega(\{x\})>0}(T^n)'(Sx)\\
&=\prod_{\omega(\{x\})>0}^{x\not\in A_S}(T^n)'(x)\prod_{\omega(\{x\})>0}^{x\in A_S}(T^n)'(Sx)\\
&=(T_*^n)'(\omega)\frac{\prod_{\omega(\{x\})>0}^{x\in A_S}(T^n)'(Sx)}
{\prod_{\omega(\{x\})>0}^{x\in A_S}(T^n)'(x)}\\
&<(T_*^n)'(\omega)\cdot  \alpha^{2\omega(A_S)}.
\endalign
$$
Hence 
$\sum_{n=1}^\infty (T_*^n)'(\omega)\le \alpha^{2\omega(A_S)}\sum_{n=1}^\infty (T_*^n)'(\omega)<+\infty$, i.e. $S_*\omega\in D$, as desired.

We now prove that the  group $\Cal S_*:=\{S_*\mid S\in\Cal S\}$ of $\mu^*$-preserving transformations of $X^*$ is ergodic\footnote{This follows also from \cite{Sh, Theorem 2.3}. However we present here an  alternative very short proof.}.
For that, we select a $\mu$-preserving totally dissipative transformation $Q$ of $X$ and a sequence $(S_n)_{n=1}^\infty$ of transformations from $\Cal S$ that converges weakly to $Q$.
Then the sequence $((S_n)_*)_{n=1}^\infty$ weakly converges to $Q_*$ in Aut$(X^*,\mu^*)$. 
If $A$ is an $\Cal S_*$-invariant subset of $X^*$ then $A$ is invariant $Q_*$ because 
$\mu^*(Q_*A\triangle A)=\lim_{n\to\infty}\mu^*((S_n)_*A\triangle A)=0$.
Since $Q_*$ is isomorphic to a probability preserving Bernoulli shift, it is ergodic and hence
$\mu^*(A)(1-\mu^*(A))=0$.
Hence $\Cal S_*$ is ergodic.

Since $D$ is $\Cal S_*$-invariant, $\mu^*(D)(1-\mu^*(D))=0$ and the assertion of the theorem follows.
\qed

\enddemo

We also prove a general result related to CCIS.

\proclaim{Proposition 7.3} Let    $T\in\text{\rom{Aut}}_1(X,\mu)$.
\roster
\item"$(i)$"
If  $(X^*,\mu^*,T_*)$ is totally dissipative then
 then $(X^*,\mu_{t}^*,T_*)$ is totally dissipative for each $t>1$.
 \item"$(ii)$"
If  $(X^*,\mu^*,T_*)$ is conservative then
 $(X^*,\mu_{t}^*,T_*)$ is conservative for each $t\in (0,1)$.
\endroster
\endproclaim

\demo{Proof} (i) If $\chi(T)\ne 0$ then the assertion of the proposition is trivial.
Hence from now on we assume that $\chi(T)=0$.
 It follows from the condition of the proposition that the Cartesian product
$
(X^*\times X^*,\mu^*\otimes\mu_{t-1}^*, T_*\times T_*).
$
is totally dissipative.
Therefore, by the Birkhoff criterion of conservativeness,
$$
\sum_{n=1}^{+\infty}\frac{d\mu^*\circ T_*^{n}}{d\mu^*}(\omega)\frac{d\mu_{t-1}^*\circ T_*^{n}}{d\mu_{t-1}^*}(\tau)<+\infty\tag 7-1
$$
at $(\mu^*\otimes\mu_{t-1}^*)$-a.e. $(\omega,\tau)\in X^*\times X^*$.
Since $T^n\in\text{\rom{Aut}}_1(X,\mu)$, it follows from \cite{DaKoRo} that for each $n>0$,
$$
\align
\frac{d\mu^*\circ T_*^{n}}{d\mu^*}(\omega)\frac{d\mu_{t-1}^*\circ T_*^{n}}{d\mu_{t-1}^*}(\tau)&=e^{\int_{X}\log \frac{d\mu\circ T^n}{d\mu}d\omega
+\int_{X}\log \frac{d\mu\circ T^n}{d\mu}d\tau}\\
&=e^{\int_{X} \log\frac{d\mu_t\circ T^n}{d\mu_t}d(\omega+\tau)}\\
&=\frac{d\mu_t^*\circ T_*^{n}}{d\mu_t^*}(\omega+\tau).
\endalign
$$
Since $\mu^**\mu_{t-1}^*=\mu^*_t$, we deduce from \thetag{7-1} that
$\sum_{n=1}^{+\infty}\frac{d\mu_t^*\circ T_*^{n}}{d\mu_t^*}(\omega)<+\infty$ at $\mu_t^*$-a.e. $\omega$.
It remains to apply the Birkhoff criterion of conservativeness.

(ii) Let $D_t:=\{\omega\in X^*\mid \sum_{n=1}^{+\infty}\frac{d\mu_t^*\circ T_*^{n}}{d\mu_t^*}(\omega)<\infty\}$ stand for the dissipative part of  $(X^*,\mu_t^*,T_*)$.
Assume in the contrapositive that $\mu_t^*(D_t)>0$.
Then the set $D_t\times X^*$ is contained in the dissipative part of the product
$(X^*\times X^*,\mu_t^*\otimes\mu_{1-t}^*,T_*\times T_*)$.
Hence for $(\mu_t\otimes\mu_{1-t})$-a.e. $(\omega,\tau)\in D_t\times X^*$,
$$
\infty>\sum_{n=1}^{+\infty}\frac{d\mu_t^*\circ T_*^{n}}{d\mu_t^*}(\omega)\frac{d\mu_{1-t}^*\circ T_*^{n}}{d\mu_{1-t}^*}(\tau)=\sum_{n=1}^{+\infty}\frac{d\mu^*\circ T_*^{n}}{d\mu^*}(\omega+\tau),
$$
i.e. $\omega+\tau$ is contained in the dissipative  part of $(X^*,\mu^*, T_*)$.
Since $\mu_t*\mu_{1-t}=\mu$, the dissipative part of $(X^*,\mu^*, T_*)$ is of  positive measure, a contradiction
\qed
\enddemo

From Proposition~7.3 and  Theorem~7.2 we deduce the main result of this section on a phase transition for conservativeness of Poisson suspensions  while   scaling the underlying  intensity.

\proclaim{Corollary 7.4} Let    $T\in\text{\rom{Aut}}_1(X,\mu)$.
Suppose that
there is  $\alpha>1$ such that $\alpha^{-1}<(T^n)'(x)<\alpha$ for each $n>0$ at a.e. $x\in X$.
Then there is $t_0\in[0,+\infty]$ such that
the Poisson suspension $(X^*,\mu_t^*,T_*)$ is conservative for each $t\in(0,t_0)$ and 
the Poisson suspension $(X^*,\mu_t^*,T_*)$ is totally dissipative for each $t\in (t_0,+\infty)$.
\endproclaim

We call $t_0$ {\it the bifurcation point}.
Of course, it is interesting when the bifurcation point is {\it proper}, i.e.
$0<t_0<+\infty$.
A concrete example of such a Poisson suspension was constructed in \S\,4 (see Remark 4.3).
Combining it with Corollary~7.4 we obtain the following theorem.

\proclaim{Theorem 7.5} Let $X=\Bbb R$, $Tx=x+1$ for all $x\in X$ and $\mu$ be an absolutely continuous measure on $X$ such that $d\mu(x)=f(x)dx$ with
$$
f(x)=\cases
1, &\text{if $x <2$}\\
-n^{-1/2},&\text{if $n\le x<n+1$ for an integer $n>1$}.
\endcases
$$
Then there exists $t_0\in[\frac16,4]$ such that
 $(X^*,\mu_t^*,T_*)$ is weakly mixing of stable type $III_1$ for each $t\in(0,t_0)$ and 
 $(X^*,\mu_t^*,T_*)$ is totally dissipative for each $t\in (t_0,+\infty)$.
\endproclaim

We note that Corollary 7.4 and Theorem 7.5 are  the ``nonsingular Poisson'' analogues of a phase transition phenomenon discovered  recently for the nonsingular Gaussian actions introduced by Y.~Arano, Y.~Isono and A.~Marrakchi in  \cite{ArIsMa} (see Proposition~5.1 there).

\head Appendix A. Skellam distributions
\endhead

{\it The Skellam distribution with parameters $(a, b)$} is the  distribution $\chi_{a,b}$ of the difference  $X-Y$ of two  independent random variables $X$ and  $Y$, each Poisson-distributed with respective parameter (expected value)  $a\ge 0$ and  $b\ge 0$. 
It is known that 
$
E(X-Y)=a-b$  and $\sigma^2(X-Y)=a+b$.
The characteristic function $\phi_{X-Y}(t):=E(e^{it(X-Y)})$ of $X-Y$ is $e^{-(a+b)+ae^{it}+be^{-it}}$ at each   $t\in\Bbb R$.
Of course, $\chi_{a,b}(\Bbb Z)=1$.
For each $k\in\Bbb Z$, 
$$
\chi_{a,b}(k)=e^{a-b}\left(\frac ab\right)^{k/2}I_k(2\sqrt{ab}),
$$
where $I_k$ is  the {\it modified Bessel function of the first kind}, i.e.
$$
I_k(z)=I_{|k|}(z):=\left(\frac z2\right)^{|k|}\sum_{j=0}^\infty\frac{(z^2/4)^j}{j!(j+|k|)!}.
$$
From the two above formulas we deduce that for each $L>0$,
$$
\align
\sum_{|k|\ge L}\chi_{a,b}(k)&\le e^{a-b}\sum_{|k|\ge L}\left(\frac ab\right)^{k/2} (ab)^{|k|/2}
\frac 1{|k|!}
\sum_{j=0}^\infty\frac{(ab)^j}{j!}\\
&=e^{a-b+ab}\left(\sum_{k\ge L}\frac{a^k}{k!}+\sum_{k\le-L}\frac{1}{b^k(-k)!}\right)\\
&=e^{a-b+ab}\left(\sum_{k\ge L}\frac{a^k}{k!}+\sum_{k\ge L}\frac{b^k}{k!}\right)\\
&=
e^{a-b+ab}\frac{a^{L}e^a+b^{L}e^b}{L!}.
\endalign
$$ 
This yields the following estimation: for each $A>0$, there is an integer $L>0$ such that for each $l>L$
$$
\sup_{0<a,b<A}\sum_{|k|\ge l}\chi_{a,b}(k)\le l^{-8}.\tag A-1
$$

\enddocument